\documentclass[11pt]{article}
\usepackage[dvips]{graphicx}
\usepackage[dvips,letterpaper,margin=1in]{geometry}
\usepackage{color}
\usepackage{amsmath}
\usepackage{amsfonts}
\usepackage{amssymb}
\usepackage{amsthm}
\usepackage{newlfont}
\usepackage{epsfig}
\usepackage{multirow}
\usepackage{setspace}
\usepackage{hyperref}
\usepackage{cite}
\usepackage{wrapfig}
\usepackage{multicol}

\begin{document}

\newtheorem{assumption}{Assumption}[section]
\newtheorem{definition}{Definition}[section]
\newtheorem{lemma}{Lemma}[section]
\newtheorem{proposition}{Proposition}[section]
\newtheorem{theorem}{Theorem}[section]
\newtheorem{corollary}{Corollary}[section]
\newtheorem{remark}{Remark}[section]
\newtheorem{conjecture}{Conjecture}[section]
\newtheorem{example}{Example}[section]

\small

\title{A note on ``MAPK networks and their capacity for multistationarity due to toric steady states''}
\author{Matthew D. Johnston\\ \\
Department of Mathematics\\University of Wisconsin-Madison\\480 Lincoln Dr., Madison, WI 53706\\email: mjohnston3@wisc.edu}
\date{}
\maketitle

\tableofcontents

\begin{abstract}
\small
We provide a short supplement to the paper ``MAPK networks and their capacity for multistationarity due to toric steady states'' by P\'{e}rez Mill\'{a}n and Turjanski. We show that the capacity for toric steady states in the three networks analyzed in that paper can be derived using the process of network translation, which corresponds the original mass action system to a generalized mass action system with the same steady states. In all three cases, the translated chemical reaction network is proper, weakly reversible, and has both a structural and kinetic deficiency of zero. This is sufficient to guarantee toric steady states by previously established work on network translations. A basis of the steady state ideal is then derived by consideration of the linkage classes of the translated chemical reaction network.
\end{abstract}


\bigskip

\section{Introduction}
\label{Introduction}

\normalsize
In \cite{P-T}, P\'{e}rez Mill\'{a}n and Turjanski consider dynamical properties of three related mitogen-activated protein kinase (MAPK) signaling pathways. They focus in particular on the capacity of the corresponding mass action systems to admit \emph{toric steady states} (Craciun \emph{et al.} \cite{C-D-S-S}; P\'{e}rez Mill\'{a}n \emph{et al.} \cite{M-D-S-C}) and on the capacity for \emph{multistationarity} and \emph{multistability} (Craciun and Feinberg \cite{C-F1,C-F2}; Feliu and Wiuf \cite{F-W2012}; M\"{u}ller and Regensberger \cite{M-R}). Multistability is the property that a mass action system admits multiple stable steady states within a single stoichiometric compatibility class. This property is of particular interest in the study of biochemical reaction networks since it is thought to underlie switching behavior. P\'{e}rez Mill\'{a}n and Turjanski use results established by P\'{e}rez Mill\'{a}n \emph{et al.} \cite{P-T} and M\"{u}ller \emph{et al.} \cite{M-F-R-C-S-D} to demonstrate that the systems admit multistability. They then numerically-derive steady states consistent with these expectations.

In this note, we elaborate on the capacity of these networks to permit toric steady states. We apply the method of network translation (Johnston \cite{J1}) to correspond the original network with a generalized network (M\"{u}ller and Regensberger \cite{M-R}) which has superior structure. In the language of \emph{chemical reaction network theory} (CRNT), we will determine that the translated network is weakly reversible and has a structural and kinetic deficiency of zero (see Feinberg \cite{F1}, Horn \cite{H}, and Horn and Jackson \cite{H-J1}; M\"{u}ller and Regensberger \cite{M-R}).

For brevity, we assume that the reader is familiar with P\'{e}rez Mill\'{a}n \emph{et al.} \cite{M-D-S-C}, M\"{u}ller and Regensberger \cite{M-R}, and Johnston \cite{J1}. We note that we use the term \emph{complex} in the sense of CRNT. That is, by a complex we mean a linear combination of reactants which appears at the terminal or root end of an arrow in the reaction graph of the network. This is distinct from the biochemical interpretation of a complex as a bound combination of two or more substrates. We adopt the following terminology and notation adapted from Johnston \cite{J1}:

\footnotesize
\begin{multicols}{2}
\begin{itemize}
\item[---]
$\mathcal{N} = (\mathcal{S},\mathcal{C},\mathcal{R})$, chemical reaction network
\item[---]
$\mathcal{N} = (\mathcal{S},\mathcal{C},\mathcal{C}_K,\mathcal{R})$, generalized chemical reaction network
\item[---]
$\tilde{\mathcal{N}} = (\tilde{\mathcal{S}},\tilde{\mathcal{C}},\tilde{\mathcal{C}}_K,\tilde{\mathcal{R}})$, translated chemical reaction network
\item[---]
$\mathcal{S} = \{ X_1, \ldots, X_m \}$, species set
\item[---]
$\mathcal{C} = \{ C_1, \ldots, C_n \}$, complex set
\item[---]
$\mathcal{C}_K = \{ (C_K)_1, \ldots, (C_K)_n \}$, kinetic complex set
\item[---]
$\mathcal{R} = \{ R_1, \ldots, R_r \}$, reaction set
\item[---]
$\mathcal{L} = \{ L_1, \ldots, L_{\ell} \}$, linkage class set
\item[---]
$m$, number of reactants
\item[---]
$n$, number of complexes
\item[---]
$r$, number of reactions
\item[---]
$\ell$, number of linkage classes
\item[---]
$\rho(i)$, index of source complex of $i^{th}$ reaction
\item[---]
$\rho'(i)$, index of product complex of $i^{th}$ reaction
\item[---]
$y_i = (y_{i1},\ldots, y_{im})$, support vector of $i^{th}$ complex
\item[---]
$(y_K)_i = ((y_K)_{i1},\ldots, (y_K)_{im})$, support vector of $i^{th}$ kinetic complex
\item[---]
$S = \mbox{span} \{ y_{\rho'(i)}-y_{\rho(i)} \; | \; i=1,\ldots, r \}$, stoichiometric subspace
\item[---]
$S_K = \mbox{span} \{ (y_K)_{\rho'(i)}-(y_K)_{\rho(i)} \; | \; i=1,\ldots, r \}$, kinetic-order subspace
\item[---]
$\mathbf{x} = (x_1, \ldots, x_m) \in \mathbb{R}_{\geq 0}^m$, concentration vector
\item[---]
$\mathbf{x}^{y_i} = \prod_{j=1}^m x_j^{y_{ij}}$, mass action term
\item[---]
$k_i,$ $i=1, \ldots, r$, reaction rate constants
\item[---]
$K_i$, $i=1,\ldots, n$, tree constants
\item[---]
$s = \mbox{dim}(S)$, dimension of stoichiometric subspace
\item[---]
$s_K = \mbox{dim}(S)$, dimension of kinetic-order subspace
\item[---]
$\delta$, deficiency ($\delta = n - \ell - s$)
\item[---]
$\delta_K$, kinetic deficiency ($\delta_K = n-\ell - s_K$)
\end{itemize}
\end{multicols}
\normalsize

\noindent We will use tilde notation $(\tilde{\cdot})$ to denote any elements above which are associated to a translated chemical reaction network. For example, we will use $\tilde{L}_i$ to denote the $i^{th}$ linkage class in the translation, $\tilde{K}_i$ to denote the tree constant corresponding to the $i^{th}$ translated complex, etc.

\section{Applications}

In this section we introduce the three networks considered by P\'{e}rez Mill\'{a}n and Turjanski \cite{P-T} and show that they admit network translations which are proper, weakly reversible, and for which both $\tilde{\delta}=0$ and $\tilde{\delta}_K = 0$. This is sufficient to guarantee toric steady states by Theorem 5 of Johnston \cite{J1}. We note that, due to similarities between the networks considered in Sections \ref{section1}, \ref{section2}, and \ref{section3} and the subsequent analysis, some clarifying exposition is omitted from Sections \ref{section2} and \ref{section3}.

\subsection{Application I: Network without feedback and three phosphatases}
\label{section1}

Consider the following network:
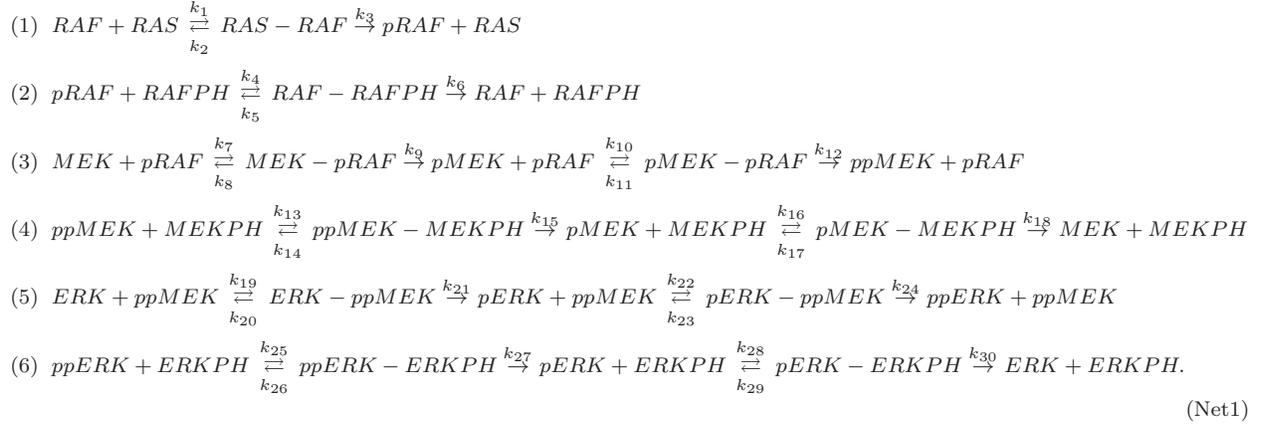
\scriptsize
\begin{equation}
\label{network1}\tag{Net1}
\begin{split}
(1) \; & \; RAF + RAS \; \mathop{\stackrel{k_1}{\rightleftarrows}}_{k_2} \; RAS-RAF \stackrel{k_3}{\rightarrow} pRAF + RAS\\
(2) \; & \; pRAF+RAFPH \; \mathop{\stackrel{k_4}{\rightleftarrows}}_{k_5} \; RAF-RAFPH \stackrel{k_6}{\rightarrow} RAF+RAFPH\\
(3) \; & \; MEK + pRAF \; \mathop{\stackrel{k_7}{\rightleftarrows}}_{k_8} \; MEK-pRAF \stackrel{k_9}{\rightarrow} pMEK+pRAF \; \mathop{\stackrel{k_{10}}{\rightleftarrows}}_{k_{11}} \; pMEK-pRAF \stackrel{k_{12}}{\rightarrow} ppMEK+pRAF\\
(4) \; & \; ppMEK+MEKPH \; \mathop{\stackrel{k_{13}}{\rightleftarrows}}_{k_{14}} \; ppMEK-MEKPH \stackrel{k_{15}}{\rightarrow} pMEK+MEKPH \; \mathop{\stackrel{k_{16}}{\rightleftarrows}}_{k_{17}} \; pMEK-MEKPH \stackrel{k_{18}}{\rightarrow} MEK + MEKPH\\
(5) \; & \; ERK+ppMEK \; \mathop{\stackrel{k_{19}}{\rightleftarrows}}_{k_{20}} \; ERK-ppMEK \stackrel{k_{21}}{\rightarrow} pERK+ppMEK \; \mathop{\stackrel{k_{22}}{\rightleftarrows}}_{k_{23}} \; pERK-ppMEK \stackrel{k_{24}}{\rightarrow} ppERK+ppMEK\\
(6) \; & \; ppERK+ERKPH \; \mathop{\stackrel{k_{25}}{\rightleftarrows}}_{k_{26}} \; ppERK-ERKPH \stackrel{k_{27}}{\rightarrow} pERK+ERKPH \; \mathop{\stackrel{k_{28}}{\rightleftarrows}}_{k_{29}} \; pERK-ERKPH \stackrel{k_{30}}{\rightarrow} ERK+ERKPH.\\
\end{split}
\end{equation}
\normalsize
The linkage classes of the network (labeled $(1)-(6)$ in (\ref{network1})) represent metabolic pathways. For instance, the pathway (1) represents the $RAS$-mediated phosphorylation of $RAF$, while pathway (2) represents the $RAFPH$-mediated dephosphorylation of $pRAF$. Note that ``$-$'' denotes a binding between the substrate on the left and the subtrate on the right. A bound substrate of this form is often called a \emph{complex} in the biochemical literature. We remind the reader, however, that we use the term complex in the sense of the CRNT literature (Horn and Jackson \cite{H-J1}).

The network (\ref{network1}) has previously been studied by Huang and Ferrell \cite{H-F1996} and Kholodenko \cite{Kholodenko2000}. The basic functionality of the mechanism may be understood by considering the pathways together as pairs: (1) \& (2), (3) \& (4), and (5) \& (6). In each such pair, there is one pathway corresponding to phorphorylation and one pathway corresponding to dephorsphorylation. In this particular model the substrates which mediate the phosphorylation and dephosphorylation events are distinct between pathways.

It was shown in P\'{e}rez Mill\'{a}n and Turjanski \cite{P-T} that the mass action system corresponding to (\ref{network1}) has toric steady states regardless of the rate constant values. Although not stated, the authors imply that the method for computing this basis was that of P\'{e}rez Mill\'{a}n \emph{et al.} \cite{M-D-S-C}. In this method, the kernel of a particular structural matrix is computed and checked for certain technical conditions (Conditions 3.1, 3.4, and 3.6, P\'{e}rez Mill\'{a}n \emph{et al.} \cite{M-D-S-C}). We now show that the basis can be derived by the method of \emph{network translation} (Johnston \cite{J1}). We note that, while the method of network translation does not produce any information in addition to that guaranteed by \cite{M-D-S-C}, it has a natural graphical representation which may be preferable to many readers.

We index the species as in P\'{e}rez Mill\'{a}n and Turjanski \cite{P-T}:

\footnotesize
\begin{equation}
\label{species1}
\begin{split}
& x_1 = [RAF], x_2 = [pRAF], x_3 = [MEK], x_4 = [pMEK], x_5 = [ppMEK], x_6 = [ERK],\\
& x_7 = [pERK], x_8 = [ppERK], x_9 = [RAS], x_{10} = [RAF PH], x_{11} = [MEKP H], x_{12} = [ERKP H],\\
& x_{13} = [RAS - RAF], x_{14} = [MEK - pRAF], x_{15} = [pMEK - pRAF], x_{16} = [ERK - ppMEK],\\
& x_{17}= [pERK - ppMEK],x_{18} = [RAF - RAF PH], x_{19} = [ppMEK - MEKP H],\\
& x_{20} = [pMEK - MEKP H], x_{21} = [ppERK - ERKP H], x_{22}= [pERK - ERKP H].
\end{split}
\end{equation}
\normalsize
We now apply the following translation scheme to (\ref{network1}) where the network has been indexed according to (\ref{species1}) and reordered by numerical index:
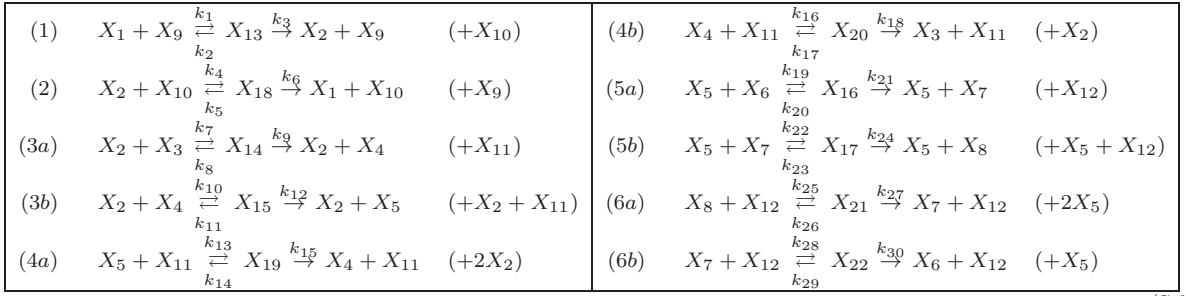
\scriptsize
\begin{equation}
\label{scheme1}\tag{Sch1}
\begin{array}{|rll|rll|}
\hline
(1) \; & \; \displaystyle{X_1 + X_9 \; \mathop{\stackrel{k_1}{\rightleftarrows}}_{k_2} \; X_{13} \stackrel{k_3}{\rightarrow} X_2 + X_9} & (+X_{10}) & 
(4b) \; & \; \displaystyle{X_4+X_{11} \; \mathop{\stackrel{k_{16}}{\rightleftarrows}}_{k_{17}} \; X_{20} \stackrel{k_{18}}{\rightarrow} X_3 + X_{11}} & (+X_{2})\\
(2) \; & \; \displaystyle{X_2+X_{10} \; \mathop{\stackrel{k_4}{\rightleftarrows}}_{k_5} \; X_{18} \stackrel{k_6}{\rightarrow} X_1+X_{10}} & (+X_{9}) &
(5a) \; & \; \displaystyle{X_5+X_6 \; \mathop{\stackrel{k_{19}}{\rightleftarrows}}_{k_{20}} \; X_{16} \stackrel{k_{21}}{\rightarrow} X_5+X_7} & (+X_{12})\\
(3a) \; & \; \displaystyle{X_2 + X_3 \; \mathop{\stackrel{k_7}{\rightleftarrows}}_{k_8} \; X_{14} \stackrel{k_9}{\rightarrow} X_2+X_4} & (+X_{11}) &
(5b) \; & \; \displaystyle{X_5+X_7 \; \mathop{\stackrel{k_{22}}{\rightleftarrows}}_{k_{23}} \; X_{17} \stackrel{k_{24}}{\rightarrow} X_5+X_8} & (+X_5 + X_{12})\\
(3b) \; & \; \displaystyle{X_2+X_4 \; \mathop{\stackrel{k_{10}}{\rightleftarrows}}_{k_{11}} \; X_{15} \stackrel{k_{12}}{\rightarrow} X_2+X_5} & (+X_{2}+X_{11}) &
(6a) \; & \; \displaystyle{X_8+X_{12} \; \mathop{\stackrel{k_{25}}{\rightleftarrows}}_{k_{26}} \; X_{21} \stackrel{k_{27}}{\rightarrow} X_7+X_{12}} & (+2X_5)\\
(4a) \; & \; \displaystyle{X_5+X_{11} \; \mathop{\stackrel{k_{13}}{\rightleftarrows}}_{k_{14}} \; X_{19} \stackrel{k_{15}}{\rightarrow} X_4+X_{11}} & (+2X_{2})&
(6b) \; & \; \displaystyle{X_7+X_{12} \; \mathop{\stackrel{k_{28}}{\rightleftarrows}}_{k_{29}} \; X_{22} \stackrel{k_{30}}{\rightarrow} X_6+X_{12}} & (+X_5)\\
\hline
\end{array}
\end{equation}
\normalsize
Notice that we have divided several individual pathways in (\ref{network1}) into separate pathways in (\ref{scheme1}). For instance, the pathway (3) has become (3a) and (3b). The purpose for these subgroupings will be apparent momentarily.

The bracketed right-hand-side terms in (\ref{scheme1}) correspond to how we translate each of the complexes in each pathway. The translation terms are chosen so that we may combine pathways to form cycles.  For example, in pathway (1) we translate the \emph{product} complex $X_1+X_9$ by $(+X_{10})$ to get the new complex $X_1 + X_9 + X_{10}$, and in pathway (2) we translate the \emph{source} complex $X_1 + X_{10}$ by $(+X_9)$ to also get $X_1 + X_9 + X_{10}$. As a result, pathway (1) and pathway (2) are joined in the translation. It is worth noting that the translation terms for the pathways (3b), (4a), (5b), and (6a) are more complicated than is strictly required to simply join the pathways together. We could, for instance, choose the terms $(+X_{11})$, $(+X_{11})$, $(+X_2)$, and $(+X_2)$ for pathways (3a), (3b), (4a), and (4b), respectively. This, however, leads to both of the \emph{source} complexes $X_2+X_4$ and $X_4 +X_{11}$ being translated to $X_2 +X_4 + X_{11}$. This overlapping of source complexes is avoided by (\ref{scheme1}) (see Application II, Johnston \cite{J1}).

The translation process (\ref{scheme1}) yields the following translated chemical reaction network $\tilde{\mathcal{N}} = (\tilde{\mathcal{S}},\tilde{\mathcal{C}},\tilde{\mathcal{C}}_K,\tilde{\mathcal{R}})$, where we have collapsed repeated complexes to a single node:

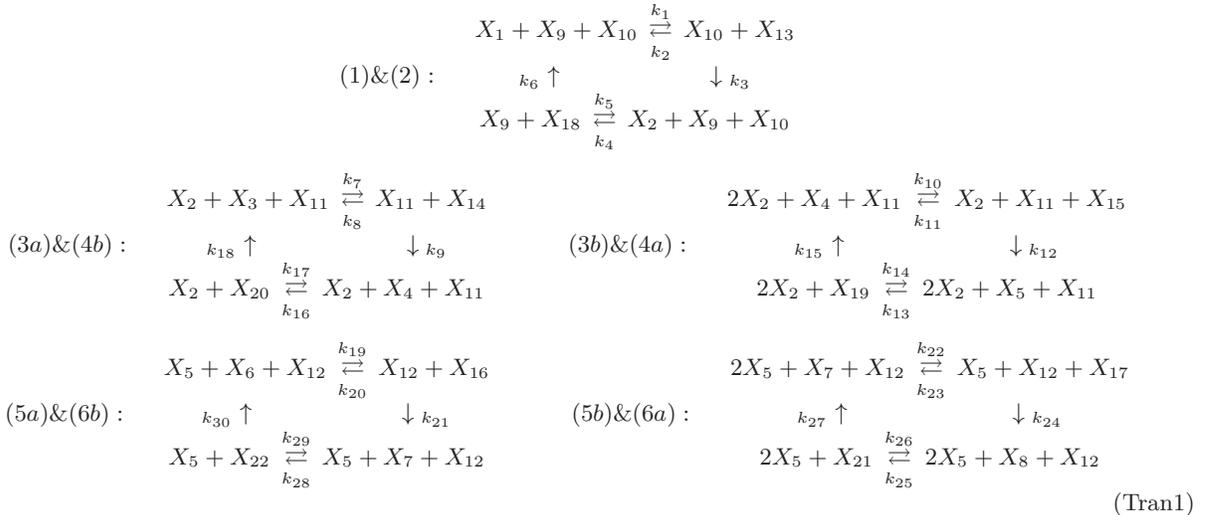
\footnotesize
\begin{equation}
\label{translation1}\tag{Tran1}
\begin{array}{c}
(1) \& (2): \; \; \; \begin{array}{c} \displaystyle{X_1 + X_9 + X_{10} \; \mathop{\stackrel{k_1}{\rightleftarrows}}_{k_2} \; X_{10} + X_{13}} \\ \displaystyle{{}_{k_6} \uparrow \hspace{2cm} \downarrow {}_{k_3}} \\ \displaystyle{X_9 + X_{18} \; \mathop{\stackrel{k_5}{\rightleftarrows}}_{k_4} \; X_2 + X_9 + X_{10}} \end{array} \vspace{0.1in} \\
(3a) \& (4b): \; \; \; \begin{array}{c} \displaystyle{X_2 + X_3 + X_{11} \; \mathop{\stackrel{k_7}{\rightleftarrows}}_{k_8} \; X_{11} + X_{14}} \\ \displaystyle{{}_{k_{18}} \uparrow \hspace{2cm} \downarrow {}_{k_{9}}} \\ \displaystyle{X_2 + X_{20} \; \mathop{\stackrel{k_{17}}{\rightleftarrows}}_{k_{16}} \; X_2 + X_4 + X_{11}} \end{array} \; \; \; \; \; \; \; \; \; \;
(3b) \& (4a): \; \; \; \begin{array}{c} \displaystyle{2X_2 + X_4 + X_{11} \; \mathop{\stackrel{k_{10}}{\rightleftarrows}}_{k_{11}} \; X_2 + X_{11} + X_{15}} \\ \displaystyle{{}_{k_{15}} \uparrow \hspace{2.2cm} \downarrow {}_{k_{12}}} \\ \displaystyle{2X_2 + X_{19} \; \mathop{\stackrel{k_{14}}{\rightleftarrows}}_{k_{13}} \; 2X_2 + X_5 + X_{11}} \end{array} \vspace{0.1in} \\
(5a) \& (6b): \; \; \; \begin{array}{c} \displaystyle{X_5 + X_6 + X_{12} \; \mathop{\stackrel{k_{19}}{\rightleftarrows}}_{k_{20}} \; X_{12} + X_{16}} \\ \displaystyle{{}_{k_{30}} \uparrow \hspace{2cm} \downarrow {}_{k_{21}}} \\ \displaystyle{X_5 + X_{22} \; \mathop{\stackrel{k_{29}}{\rightleftarrows}}_{k_{28}} \; X_5 + X_7 + X_{12}} \end{array} \; \; \; \; \; \; \; \; \; \;
(5b) \& (6a): \; \; \; \begin{array}{c} \displaystyle{2X_5 + X_7 + X_{12} \; \mathop{\stackrel{k_{22}}{\rightleftarrows}}_{k_{23}} \; X_5 + X_{12} + X_{17}} \\ \displaystyle{{}_{k_{27}} \uparrow \hspace{2.2cm} \downarrow {}_{k_{24}}} \\ \displaystyle{2X_5 + X_{21} \; \mathop{\stackrel{k_{26}}{\rightleftarrows}}_{k_{25}} \; 2X_5 + X_8 + X_{12}} \end{array}
\end{array}
\end{equation}
\normalsize
The kinetic complex set $\tilde{\mathcal{C}}_K$ is omitted from (\ref{translation1}) to conserve space. We now analyze (\ref{translation1}) according to the methodologies introduced by Johnston \cite{J1}. We first construct the kinetic graph $(\tilde{\mathcal{S}},\tilde{\mathcal{C}}_K,\tilde{\mathcal{R}})$:

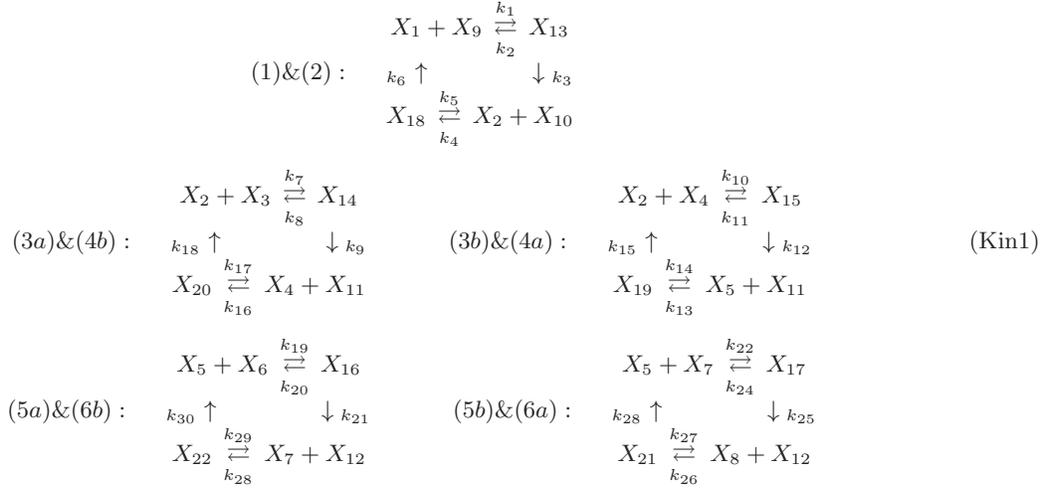
\footnotesize
\begin{equation}
\label{kinetic1}\tag{Kin1}
\begin{array}{c}
(1) \& (2): \; \; \; \begin{array}{c} \displaystyle{X_1 + X_9 \; \mathop{\stackrel{k_1}{\rightleftarrows}}_{k_2} \; X_{13}} \\ \displaystyle{{}_{k_6} \uparrow \hspace{1.4cm} \downarrow {}_{k_3}} \\ \displaystyle{X_{18} \; \mathop{\stackrel{k_5}{\rightleftarrows}}_{k_4} \; X_2 + X_{10}} \end{array} \vspace{0.1in} \\
(3a) \& (4b): \; \; \; \begin{array}{c} \displaystyle{X_2 + X_3 \; \mathop{\stackrel{k_7}{\rightleftarrows}}_{k_8} \; X_{14}} \\ \displaystyle{{}_{k_{18}} \uparrow \hspace{1.4cm} \downarrow {}_{k_{9}}} \\ \displaystyle{X_{20} \; \mathop{\stackrel{k_{17}}{\rightleftarrows}}_{k_{16}} \; X_4 + X_{11}} \end{array} \; \; \; \; \; \; \; \; \; \;
(3b) \& (4a): \; \; \; \begin{array}{c} \displaystyle{X_2 + X_4 \; \mathop{\stackrel{k_{10}}{\rightleftarrows}}_{k_{11}} \; X_{15}} \\ \displaystyle{{}_{k_{15}} \uparrow \hspace{1.4cm} \downarrow {}_{k_{12}}} \\ \displaystyle{X_{19} \; \mathop{\stackrel{k_{14}}{\rightleftarrows}}_{k_{13}} \; X_5 + X_{11}} \end{array} \vspace{0.1in} \\
(5a) \& (6b): \; \; \; \begin{array}{c} \displaystyle{X_5 + X_6 \; \mathop{\stackrel{k_{19}}{\rightleftarrows}}_{k_{20}} \; X_{16}} \\ \displaystyle{{}_{k_{30}} \uparrow \hspace{1.4cm} \downarrow {}_{k_{21}}} \\ \displaystyle{X_{22} \; \mathop{\stackrel{k_{29}}{\rightleftarrows}}_{k_{28}} \; X_7 + X_{12}} \end{array} \; \; \; \; \; \; \; \; \; \;
(5b) \& (6a): \; \; \; \begin{array}{c} \displaystyle{X_5 + X_7 \; \mathop{\stackrel{k_{22}}{\rightleftarrows}}_{k_{24}} \; X_{17}} \\ \displaystyle{{}_{k_{28}} \uparrow \hspace{1.4cm} \downarrow {}_{k_{25}}} \\ \displaystyle{X_{21} \; \mathop{\stackrel{k_{27}}{\rightleftarrows}}_{k_{26}} \; X_8 + X_{12}} \end{array}
\end{array}
\end{equation}
\normalsize
The network (\ref{kinetic1}) is determined by substituting the pre-translation source complexes in (\ref{scheme1}) in the place of the corresponding post-translation source complexes in (\ref{translation1}). Since the correspondence between these source complexes is one-to-one, the translation is proper so the dynamics of the original mass action system corresponding to (\ref{network1}) coincides with that of the generalized system corresponding to the network (\ref{translation1}) with kinetic complexes given by (\ref{kinetic1}) (Lemma 2, Johnston \cite{J1}). The generalized network $\tilde{\mathcal{N}}$ is weakly reversible and has $\tilde{\delta}=0$ (zero deficiency of (\ref{translation1})) and $\tilde{\delta}_K = 0$ (zero deficiency of (\ref{kinetic1})). It follows that the system has toric steady states regardless of rate constant values (Property 1. of Theorem 5, Johnston \cite{J1}).

An explicit derivation of the toric basis of the steady state ideal given by P\'{e}rez Mill\'{a}n and Turjanski \cite{P-T} can now be made by consideration of the kinetic graph (\ref{kinetic1}). By property 2. of Theorem 5 of Johnston \cite{J1}, the terms we require are
\[\tilde{K}_i \mathbf{x}^{(\tilde{y}_K)_j} - \tilde{K}_j \mathbf{x}^{(\tilde{y}_K)_i}\]
where $i, j \in \tilde{L}_k$ for some linkage class $\tilde{L}_k$ of (\ref{kinetic1}), and $\tilde{K}_i$ are the tree constants of the translation. 
We make the following enumeration of the kinetic complexes:

\footnotesize
\begin{equation}
\label{kineticcomplexes1}
\begin{split}
(1) \& (2): \; \; \; \; & \; (\tilde{C}_K)_1 = X_1 + X_9, (\tilde{C}_K)_2 = X_{13}, (\tilde{C}_K)_3 = X_2+X_{10}, (\tilde{C}_K)_4 = X_{18},\\
(3a) \& (4b): \; \; \; \; & \; (\tilde{C}_K)_5 = X_2 + X_3, (\tilde{C}_K)_6 = X_{14}, (\tilde{C}_K)_7 = X_4 + X_{11}, (\tilde{C}_K)_8 = X_{20},\\
(3b) \& (4a): \; \; \; \; & \; (\tilde{C}_K)_9 = X_2 + X_4, (\tilde{C}_K)_{10} = X_{15}, (\tilde{C}_K)_{11} = X_5 + X_{11}, (\tilde{C}_K)_{12} = X_{19},\\
(5a) \& (6b): \; \; \; \; & \; (\tilde{C}_K)_{13} = X_5 + X_6, (\tilde{C}_K)_{14} = X_{16}, (\tilde{C}_K)_{15} = X_7 + X_{12}, (\tilde{C}_K)_{16} = X_{22},\\
(5b) \& (6a): \; \; \; \; & \; (\tilde{C}_K)_{17} = X_5 + X_7, (\tilde{C}_K)_{18} = X_{17}, (\tilde{C}_K)_{19} = X_8 + X_{12}, (\tilde{C}_K)_{20} = X_{21}
\end{split}
\end{equation}
\normalsize
where the complexes have been grouped according to their linkage class in (\ref{kinetic1}). We compute the following tree constants:

\footnotesize
\begin{equation}
\begin{array}{|l|l|l|l|}
\hline
\tilde{K}_1 = (k_2 + k_3)k_4k_6 & \tilde{K}_6 = k_7k_{16}k_{18} & \tilde{K}_{11} = k_{10}k_{12}(k_{14}+k_{15}) & \tilde{K}_{16} = k_{19}k_{21}k_{28}\\
\tilde{K}_2 = k_1k_4k_6 & \tilde{K}_7 = k_7k_9(k_{17}+k_{18}) & \tilde{K}_{12} = k_{10}k_{12}k_{13} & \tilde{K}_{17} = (k_{23}+k_{24})k_{25}k_{27}\\
\tilde{K}_3 = k_1k_3(k_5+k_6) & \tilde{K}_8 = k_7k_9k_{16} & \tilde{K}_{13} = (k_{20}+k_{21})k_{28}k_{30} & \tilde{K}_{18} = k_{22}k_{25}k_{27}\\
\tilde{K}_4 = k_1k_3k_4 & \tilde{K}_9 = (k_{11}+k_{12})k_{13}k_{15} & \tilde{K}_{14} = k_{19}k_{28}k_{30} & \tilde{K}_{19} = k_{22}k_{24}(k_{26}+k_{27})\\
\tilde{K}_5 = (k_8+k_9)k_{16}k_{18} & \tilde{K}_{10} = k_{10}k_{13}k_{15} &\tilde{K}_{15} = k_{19}k_{21}(k_{29}+k_{30}) & \tilde{K}_{20} = k_{23}k_{24}k_{25}\\
\hline
\end{array}
\end{equation}
\normalsize
In accordance with the results of P\'{e}rez Mill\'{a}n and Turjanski \cite{P-T}, after simplification we obtain the following:
\footnotesize
\[\begin{array}{c}
\begin{array}{c}
\underline{(1) \& (2)}: \vspace{0.05in} \\
\displaystyle{\tilde{K}_1\mathbf{x}^{(\tilde{y}_K)_2} - \tilde{K}_2 \mathbf{x}^{(\tilde{y}_K)_1} \Rightarrow (k_2 + k_3) x_{13} - k_1 x_1x_9}\\
\displaystyle{\tilde{K}_3\mathbf{x}^{(\tilde{y}_K)_4} - \tilde{K}_4 \mathbf{x}^{(\tilde{y}_K)_3} \Rightarrow (k_5 + k_6) x_{18} - k_4 x_2x_{10}}\\
\displaystyle{\tilde{K}_4\mathbf{x}^{(\tilde{y}_K)_2} - \tilde{K}_2 \mathbf{x}^{(\tilde{y}_K)_4} \Rightarrow k_3 x_{13} - k_6 x_{18}}
\end{array} \vspace{0.1in}\\
\begin{array}{c}
\underline{(3a) \& (4b)}: \vspace{0.05in} \\
\displaystyle{\tilde{K}_5\mathbf{x}^{(\tilde{y}_K)_6} - \tilde{K}_6 \mathbf{x}^{(\tilde{y}_K)_5} \Rightarrow (k_8 + k_9) x_{14} - k_7 x_2x_3}\\
\displaystyle{\tilde{K}_7\mathbf{x}^{(\tilde{y}_K)_8} - \tilde{K}_8 \mathbf{x}^{(\tilde{y}_K)_7} \Rightarrow (k_{17} + k_{18}) x_{20} - k_{16} x_4x_{11}}\\
\displaystyle{\tilde{K}_8\mathbf{x}^{(\tilde{y}_K)_6} - \tilde{K}_6 \mathbf{x}^{(\tilde{y}_K)_8} \Rightarrow k_{9} x_{14} - k_{18} x_{20}}
\end{array} \hspace{0.1in}
\begin{array}{c}
\underline{(3b) \& (4a)}: \vspace{0.05in} \\
\displaystyle{\tilde{K}_9\mathbf{x}^{(\tilde{y}_K)_{10}}- \tilde{K}_{10} \mathbf{x}^{(\tilde{y}_K)_9} \Rightarrow (k_{11} + k_{12}) x_{15} - k_{10} x_2x_4}\\
\displaystyle{\tilde{K}_{11}\mathbf{x}^{(\tilde{y}_K)_{12}} - \tilde{K}_{12} \mathbf{x}^{(\tilde{y}_K)_{11}} \Rightarrow (k_{14} + k_{15}) x_{19} - k_{13} x_5x_{11}}\\
\displaystyle{\tilde{K}_{12}\mathbf{x}^{(\tilde{y}_K)_{10}} - \tilde{K}_{10} \mathbf{x}^{(\tilde{y}_K)_{12}} \Rightarrow k_{12} x_{15} - k_{15} x_{19}}
\end{array} \vspace{0.1in}\\
\begin{array}{c}
\underline{(5a) \& (6b)}: \vspace{0.05in} \\
\displaystyle{\tilde{K}_{13}\mathbf{x}^{(\tilde{y}_K)_{14}} - \tilde{K}_{14} \mathbf{x}^{(\tilde{y}_K)_{13}} \Rightarrow (k_{20} + k_{21}) x_{16} - k_{19} x_5x_6}\\
\displaystyle{\tilde{K}_{15}\mathbf{x}^{(\tilde{y}_K)_{16}} - \tilde{K}_{16} \mathbf{x}^{(\tilde{y}_K)_{15}} \Rightarrow (k_{29} + k_{30}) x_{22} - k_{28} x_2x_{12}}\\
\displaystyle{\tilde{K}_{16}\mathbf{x}^{(\tilde{y}_K)_{14}} - \tilde{K}_{14} \mathbf{x}^{(\tilde{y}_K)_{16}} \Rightarrow k_{21} x_{16} - k_{30} x_{22}}
\end{array} \hspace{0.1in}
\begin{array}{c}
\underline{(5b) \& (6a)}: \vspace{0.05in} \\
\displaystyle{\tilde{K}_{17}\mathbf{x}^{(\tilde{y}_K)_{18}} - \tilde{K}_{18} \mathbf{x}^{(\tilde{y}_K)_{17}} \Rightarrow (k_{23} + k_{24}) x_{17} - k_{22} x_5x_7}\\
\displaystyle{\tilde{K}_{19}\mathbf{x}^{(\tilde{y}_K)_{20}} - \tilde{K}_{20} \mathbf{x}^{(\tilde{y}_K)_{19}} \Rightarrow (k_{26} + k_{27}) x_{21} - k_{25} x_8x_{12}}\\
\displaystyle{\tilde{K}_{20}\mathbf{x}^{(\tilde{y}_K)_{20}} - \tilde{K}_{18} \mathbf{x}^{(\tilde{y}_K)_{20}} \Rightarrow k_{24} x_{17} - k_{27} x_{21}}
\end{array}
\end{array}\]
\normalsize

The elements of this basis are grouped to emphasize the connection with the linkage classes of the translation (\ref{translation1}). The basis can be constructed by considering any independent pair-wise combination of complexes in a given linkage class of (\ref{kinetic1}), computing, and simplifying $\tilde{K}_i \mathbf{x}^{(\tilde{y}_K)_j} - \tilde{K}_j \mathbf{x}^{(\tilde{y}_K)_i}$. In general, we have $\tilde{n}-\tilde{\ell}$ such choices where $\tilde{n}$ is the number of complexes in the translation and $\tilde{\ell}$ is the number of linkage classes. For this network, we need three basis pairs from each linkage class, since there are four complexes in each linkage class, which gives a total of $15$ independent basis pairs since there are five linkage classes.

\subsection{Application II: Network without feedback and two phosphatases}
\label{section2}

Consider the following network:
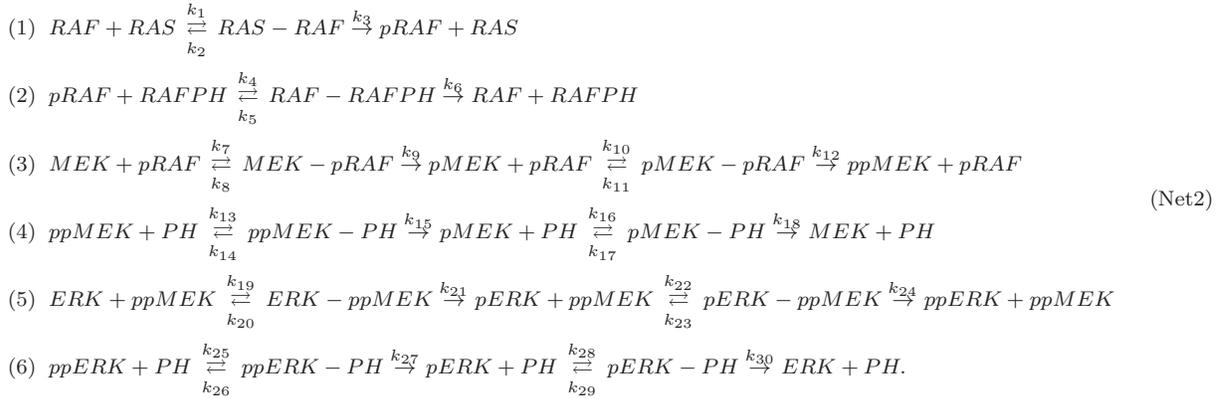
\scriptsize
\begin{equation}
\label{network2}\tag{Net2}
\begin{split}
(1) \; & \; RAF + RAS \; \mathop{\stackrel{k_1}{\rightleftarrows}}_{k_2} \; RAS-RAF \stackrel{k_3}{\rightarrow} pRAF + RAS\\
(2) \; & \; pRAF+RAFPH \; \mathop{\stackrel{k_4}{\rightleftarrows}}_{k_5} \; RAF-RAFPH \stackrel{k_6}{\rightarrow} RAF+RAFPH\\
(3) \; & \; MEK + pRAF \; \mathop{\stackrel{k_7}{\rightleftarrows}}_{k_8} \; MEK-pRAF \stackrel{k_9}{\rightarrow} pMEK+pRAF \; \mathop{\stackrel{k_{10}}{\rightleftarrows}}_{k_{11}} \; pMEK-pRAF \stackrel{k_{12}}{\rightarrow} ppMEK+pRAF\\
(4) \; & \; ppMEK+PH \; \mathop{\stackrel{k_{13}}{\rightleftarrows}}_{k_{14}} \; ppMEK-PH \stackrel{k_{15}}{\rightarrow} pMEK+PH \; \mathop{\stackrel{k_{16}}{\rightleftarrows}}_{k_{17}} \; pMEK-PH \stackrel{k_{18}}{\rightarrow} MEK + PH\\
(5) \; & \; ERK+ppMEK \; \mathop{\stackrel{k_{19}}{\rightleftarrows}}_{k_{20}} \; ERK-ppMEK \stackrel{k_{21}}{\rightarrow} pERK+ppMEK \; \mathop{\stackrel{k_{22}}{\rightleftarrows}}_{k_{23}} \; pERK-ppMEK \stackrel{k_{24}}{\rightarrow} ppERK+ppMEK\\
(6) \; & \; ppERK+PH \; \mathop{\stackrel{k_{25}}{\rightleftarrows}}_{k_{26}} \; ppERK-PH \stackrel{k_{27}}{\rightarrow} pERK+PH \; \mathop{\stackrel{k_{28}}{\rightleftarrows}}_{k_{29}} \; pERK-PH \stackrel{k_{30}}{\rightarrow} ERK+PH.\\
\end{split}
\end{equation}
\normalsize
The network is identical to (\ref{network1}) except that the dephosphorylation pathways (4) and (6) are mediated by the same phosphotase (denoted $PH$). This mechanism has been previously considered by Fujioka \emph{et al.} \cite{Fujioka2006}.

We now perform the same analysis as contained in Section \ref{section1}. We start by indexing the reactants as in P\'{e}rez Mill\'{a}n and Turjanski \cite{P-T}:

\footnotesize
\begin{equation}
\label{species2}
\begin{split}
& x_1 = [RAF], x_2 = [pRAF], x_3 = [MEK], x_4 = [pMEK], x_5 = [ppMEK], x_6 = [ERK], x_7 = [pERK],\\
& x_8 = [ppERK], x_9 = [RAS], x_{10} = [RAF P H], x_{11} = [P H], x_{12} = [RAS - RAF], x_{13} = [MEK - pRAF],\\
& x_{14} = [pMEK - pRAF], x_{15} = [ERK - ppMEK], x_{16} = [pERK - ppMEK], x_{17} = [RAF - RAF PH],\\
& x_{18} = [ppMEK - P H], x_{19} = [pMEK - P H], x_{20} = [ppERK - P H], x_{21} = [pERK - P H].
\end{split}
\end{equation}
\normalsize
We now apply the following translation scheme where the reactants have been indexed according to (\ref{species2}):

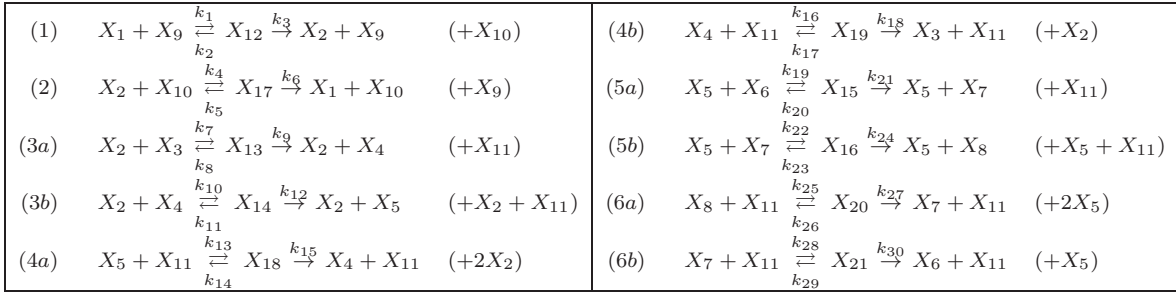
\scriptsize
\begin{equation}
\label{scheme2}\tag{Sch2}
\begin{array}{|rll|rll|}
\hline
(1) \; & \; \displaystyle{X_1 + X_9 \; \mathop{\stackrel{k_1}{\rightleftarrows}}_{k_2} \; X_{12} \stackrel{k_3}{\rightarrow} X_2 + X_9} & (+X_{10}) & 
(4b) \; & \; \displaystyle{X_4+X_{11} \; \mathop{\stackrel{k_{16}}{\rightleftarrows}}_{k_{17}} \; X_{19} \stackrel{k_{18}}{\rightarrow} X_3 + X_{11}} & (+X_{2})\\
(2) \; & \; \displaystyle{X_2+X_{10} \; \mathop{\stackrel{k_4}{\rightleftarrows}}_{k_5} \; X_{17} \stackrel{k_6}{\rightarrow} X_1+X_{10}} & (+X_{9}) &
(5a) \; & \; \displaystyle{X_5+X_6 \; \mathop{\stackrel{k_{19}}{\rightleftarrows}}_{k_{20}} \; X_{15} \stackrel{k_{21}}{\rightarrow} X_5+X_7} & (+X_{11})\\
(3a) \; & \; \displaystyle{X_2 + X_3 \; \mathop{\stackrel{k_7}{\rightleftarrows}}_{k_8} \; X_{13} \stackrel{k_9}{\rightarrow} X_2+X_4} & (+X_{11}) &
(5b) \; & \; \displaystyle{X_5+X_7 \; \mathop{\stackrel{k_{22}}{\rightleftarrows}}_{k_{23}} \; X_{16} \stackrel{k_{24}}{\rightarrow} X_5+X_8} & (+X_5 + X_{11})\\
(3b) \; & \; \displaystyle{X_2+X_4 \; \mathop{\stackrel{k_{10}}{\rightleftarrows}}_{k_{11}} \; X_{14} \stackrel{k_{12}}{\rightarrow} X_2+X_5} & (+X_{2}+X_{11}) &
(6a) \; & \; \displaystyle{X_8+X_{11} \; \mathop{\stackrel{k_{25}}{\rightleftarrows}}_{k_{26}} \; X_{20} \stackrel{k_{27}}{\rightarrow} X_7+X_{11}} & (+2X_5)\\
(4a) \; & \; \displaystyle{X_5+X_{11} \; \mathop{\stackrel{k_{13}}{\rightleftarrows}}_{k_{14}} \; X_{18} \stackrel{k_{15}}{\rightarrow} X_4+X_{11}} & (+2X_{2})&
(6b) \; & \; \displaystyle{X_7+X_{11} \; \mathop{\stackrel{k_{28}}{\rightleftarrows}}_{k_{29}} \; X_{21} \stackrel{k_{30}}{\rightarrow} X_6+X_{11}} & (+X_5)\\
\hline
\end{array}
\end{equation}
\normalsize
The translation process (\ref{scheme2}) yields the following translated chemical reaction network $\tilde{\mathcal{N}} = (\tilde{\mathcal{S}},\tilde{\mathcal{C}},\tilde{\mathcal{C}}_K,\tilde{\mathcal{R}})$:
\footnotesize
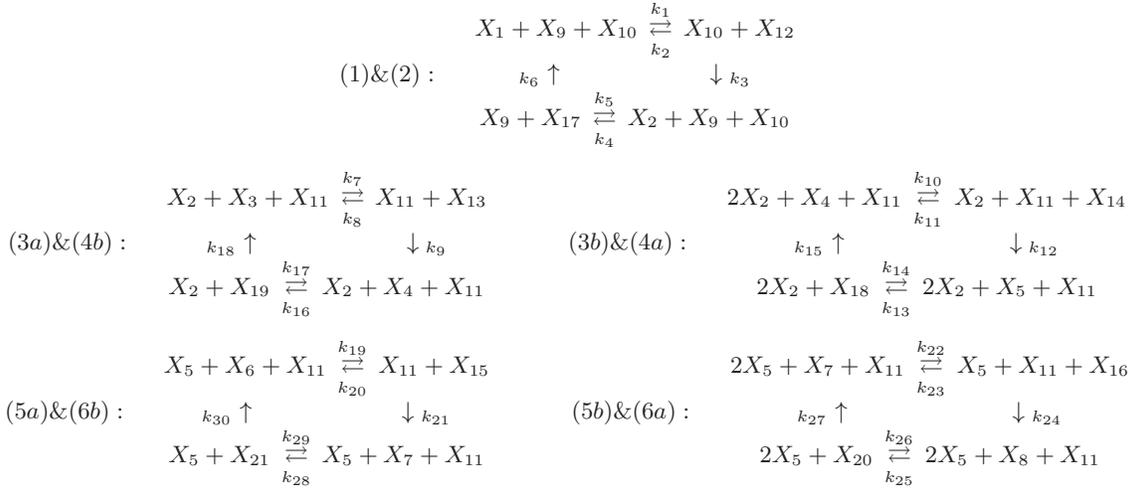
\begin{equation}
\label{translation2}\tag{Tran2}
\begin{array}{c}
(1) \& (2): \; \; \; \begin{array}{c} \displaystyle{X_1 + X_9 + X_{10} \; \mathop{\stackrel{k_1}{\rightleftarrows}}_{k_2} \; X_{10} + X_{12}} \\ \displaystyle{{}_{k_6} \uparrow \hspace{2cm} \downarrow {}_{k_3}} \\ \displaystyle{X_9 + X_{17} \; \mathop{\stackrel{k_5}{\rightleftarrows}}_{k_4} \; X_2 + X_9 + X_{10}} \end{array} \vspace{0.1in} \\
(3a) \& (4b): \; \; \; \begin{array}{c} \displaystyle{X_2 + X_3 + X_{11} \; \mathop{\stackrel{k_7}{\rightleftarrows}}_{k_8} \; X_{11} + X_{13}} \\ \displaystyle{{}_{k_{18}} \uparrow \hspace{2cm} \downarrow {}_{k_{9}}} \\ \displaystyle{X_2 + X_{19} \; \mathop{\stackrel{k_{17}}{\rightleftarrows}}_{k_{16}} \; X_2 + X_4 + X_{11}} \end{array} \; \; \; \; \; \; \; \; \; \;
(3b) \& (4a): \; \; \; \begin{array}{c} \displaystyle{2X_2 + X_4 + X_{11} \; \mathop{\stackrel{k_{10}}{\rightleftarrows}}_{k_{11}} \; X_2 + X_{11} + X_{14}} \\ \displaystyle{{}_{k_{15}} \uparrow \hspace{2.2cm} \downarrow {}_{k_{12}}} \\ \displaystyle{2X_2 + X_{18} \; \mathop{\stackrel{k_{14}}{\rightleftarrows}}_{k_{13}} \; 2X_2 + X_5 + X_{11}} \end{array} \vspace{0.1in} \\
(5a) \& (6b): \; \; \; \begin{array}{c} \displaystyle{X_5 + X_6 + X_{11} \; \mathop{\stackrel{k_{19}}{\rightleftarrows}}_{k_{20}} \; X_{11} + X_{15}} \\ \displaystyle{{}_{k_{30}} \uparrow \hspace{2cm} \downarrow {}_{k_{21}}} \\ \displaystyle{X_5 + X_{21} \; \mathop{\stackrel{k_{29}}{\rightleftarrows}}_{k_{28}} \; X_5 + X_7 + X_{11}} \end{array} \; \; \; \; \; \; \; \; \; \;
(5b) \& (6a): \; \; \; \begin{array}{c} \displaystyle{2X_5 + X_7 + X_{11} \; \mathop{\stackrel{k_{22}}{\rightleftarrows}}_{k_{23}} \; X_5 + X_{11} + X_{16}} \\ \displaystyle{{}_{k_{27}} \uparrow \hspace{2.2cm} \downarrow {}_{k_{24}}} \\ \displaystyle{2X_5 + X_{20} \; \mathop{\stackrel{k_{26}}{\rightleftarrows}}_{k_{25}} \; 2X_5 + X_8 + X_{11}} \end{array}
\end{array}
\end{equation}
\normalsize
and the associated kinetic graph $(\tilde{\mathcal{S}},\tilde{\mathcal{C}}_K,\tilde{\mathcal{R}})$:
\footnotesize
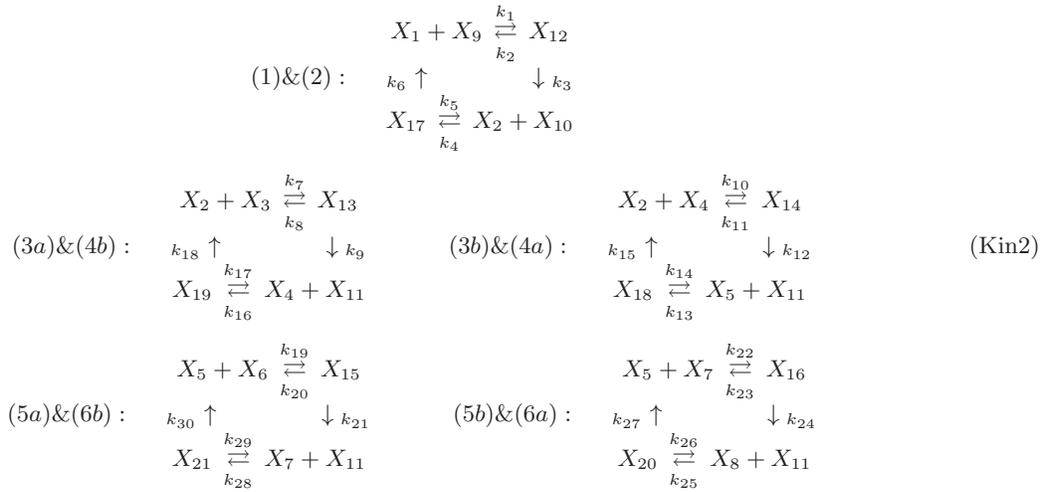
\begin{equation}
\label{kinetic2}\tag{Kin2}
\begin{array}{c}
(1) \& (2): \; \; \; \begin{array}{c} \displaystyle{X_1 + X_9 \; \mathop{\stackrel{k_1}{\rightleftarrows}}_{k_2} \; X_{12}} \\ \displaystyle{{}_{k_6} \uparrow \hspace{1.4cm} \downarrow {}_{k_3}} \\ \displaystyle{X_{17} \; \mathop{\stackrel{k_5}{\rightleftarrows}}_{k_4} \; X_2 + X_{10}} \end{array} \vspace{0.1in} \\
(3a) \& (4b): \; \; \; \begin{array}{c} \displaystyle{X_2 + X_3 \; \mathop{\stackrel{k_7}{\rightleftarrows}}_{k_8} \; X_{13}} \\ \displaystyle{{}_{k_{18}} \uparrow \hspace{1.4cm} \downarrow {}_{k_{9}}} \\ \displaystyle{X_{19} \; \mathop{\stackrel{k_{17}}{\rightleftarrows}}_{k_{16}} \; X_4 + X_{11}} \end{array} \; \; \; \; \; \; \; \; \; \;
(3b) \& (4a): \; \; \; \begin{array}{c} \displaystyle{X_2 + X_4 \; \mathop{\stackrel{k_{10}}{\rightleftarrows}}_{k_{11}} \; X_{14}} \\ \displaystyle{{}_{k_{15}} \uparrow \hspace{1.4cm} \downarrow {}_{k_{12}}} \\ \displaystyle{X_{18} \; \mathop{\stackrel{k_{14}}{\rightleftarrows}}_{k_{13}} \; X_5 + X_{11}} \end{array} \vspace{0.1in} \\
(5a) \& (6b): \; \; \; \begin{array}{c} \displaystyle{X_5 + X_6 \; \mathop{\stackrel{k_{19}}{\rightleftarrows}}_{k_{20}} \; X_{15}} \\ \displaystyle{{}_{k_{30}} \uparrow \hspace{1.4cm} \downarrow {}_{k_{21}}} \\ \displaystyle{X_{21} \; \mathop{\stackrel{k_{29}}{\rightleftarrows}}_{k_{28}} \; X_7 + X_{11}} \end{array} \; \; \; \; \; \; \; \; \; \;
(5b) \& (6a): \; \; \; \begin{array}{c} \displaystyle{X_5 + X_7 \; \mathop{\stackrel{k_{22}}{\rightleftarrows}}_{k_{23}} \; X_{16}} \\ \displaystyle{{}_{k_{27}} \uparrow \hspace{1.4cm} \downarrow {}_{k_{24}}} \\ \displaystyle{X_{20} \; \mathop{\stackrel{k_{26}}{\rightleftarrows}}_{k_{25}} \; X_8 + X_{11}} \end{array}
\end{array}
\end{equation}
\normalsize
As for the network in Section \ref{section1}, we can determine that the translation $\tilde{\mathcal{N}}$ is proper, weakly reversible, and has both $\tilde{\delta} = 0$ and $\tilde{\delta}_K = 0$. We make the following enumeration of the kinetic complexes:

\footnotesize
\begin{equation}
\label{kineticcomplexes2}
\begin{split}
(1) \& (2): \; \; \; \; \; & (\tilde{C}_K)_1 = X_1 + X_9, (\tilde{C}_K)_2 = X_{12}, (\tilde{C}_K)_3 = X_2+X_{10}, (\tilde{C}_K)_4 = X_{17},\\
(3a) \& (4b): \; \; \; \; \; & (\tilde{C}_K)_5 = X_2 + X_3, (\tilde{C}_K)_6 = X_{13}, (\tilde{C}_K)_7 = X_4 + X_{11}, (\tilde{C}_K)_8 = X_{19},\\
(3b) \& (4a): \; \; \; \; \; & (\tilde{C}_K)_9 = X_2 + X_4, (\tilde{C}_K)_{10} = X_{14}, (\tilde{C}_K)_{11} = X_5 + X_{11}, (\tilde{C}_K)_{12} = X_{18},\\
(5a) \& (6b): \; \; \; \; \; & (\tilde{C}_K)_{13} = X_5 + X_6, (\tilde{C}_K)_{14} = X_{15}, (\tilde{C}_K)_{15} = X_7 + X_{11}, (\tilde{C}_K)_{16} = X_{21},\\
(5b) \& (6a): \; \; \; \; \; & (\tilde{C}_K)_{17} = X_5 + X_7, (\tilde{C}_K)_{18} = X_{16}, (\tilde{C}_K)_{19} = X_8 + X_{11}, (\tilde{C}_K)_{20} = X_{20},
\end{split}
\end{equation}
\normalsize
and compute the following tree constants:

\footnotesize
\begin{equation}
\begin{array}{|l|l|l|l|}
\hline
\tilde{K}_1 = (k_2 + k_3)k_4k_6 & \tilde{K}_6 = k_7k_{18}k_{16} & \tilde{K}_{11} = k_{10}k_{12}(k_{14}+k_{15}) & \tilde{K}_{16} = k_{19}k_{21}k_{28}\\
\tilde{K}_2 = k_1k_4k_6 & \tilde{K}_7 = k_7k_9(k_{17}+k_{18}) & \tilde{K}_{12} = k_{10}k_{12}k_{13} & \tilde{K}_{17} = (k_{23}+k_{24})k_{25}k_{28}\\
\tilde{K}_3 = k_1k_3(k_5+k_6) & \tilde{K}_8 = k_7k_9k_{16} & \tilde{K}_{13} = (k_{20}+k_{21})k_{28}k_{30} & \tilde{K}_{18} = k_{22}k_{25}k_{27}\\
\tilde{K}_4 = k_1k_3k_4 & \tilde{K}_9 = (k_{11}+k_{12})k_{13}k_{15} & \tilde{K}_{14} = k_{19}k_{28}k_{30} & \tilde{K}_{19} = k_{22}k_{24}(k_{26}+k_{27})\\
\tilde{K}_5 = (k_8+k_9)k_{16}k_{18} & \tilde{K}_{10} = k_{10}k_{13}k_{15} &\tilde{K}_{15} = k_{19}k_{21}(k_{29}+k_{30}) & \tilde{K}_{20} = k_{23}k_{24}k_{25}\\
\hline
\end{array}
\end{equation}
\normalsize
By property 2. of Theorem 5 of Johnston \cite{J1}, a basis of the steady state ideal can then be given by:

\footnotesize
\[\begin{array}{c}
\begin{array}{c}
\underline{(1) \& (2)}: \vspace{0.05in} \\
\displaystyle{\tilde{K}_1\mathbf{x}^{(\tilde{y}_K)_2} - \tilde{K}_2 \mathbf{x}^{(\tilde{y}_K)_1} \Rightarrow (k_2 + k_3) x_{12} - k_1 x_1x_9}\\
\displaystyle{\tilde{K}_3\mathbf{x}^{(\tilde{y}_K)_4} - \tilde{K}_4 \mathbf{x}^{(\tilde{y}_K)_3} \Rightarrow (k_5 + k_6) x_{17} - k_4 x_2x_{10}}\\
\displaystyle{\tilde{K}_4\mathbf{x}^{(\tilde{y}_K)_2} - \tilde{K}_2 \mathbf{x}^{(\tilde{y}_K)_4} \Rightarrow k_3 x_{12} - k_6 x_{17}}
\end{array} \vspace{0.1in}\\
\begin{array}{c}
\underline{(3a) \& (4b)}: \vspace{0.05in} \\
\displaystyle{\tilde{K}_5\mathbf{x}^{(\tilde{y}_K)_6} - \tilde{K}_6 \mathbf{x}^{(\tilde{y}_K)_5} \Rightarrow (k_8 + k_9) x_{13} - k_7 x_2x_3}\\
\displaystyle{\tilde{K}_7\mathbf{x}^{(\tilde{y}_K)_8} - \tilde{K}_8 \mathbf{x}^{(\tilde{y}_K)_7} \Rightarrow (k_{17} + k_{18}) x_{19} - k_{16} x_4x_{11}}\\
\displaystyle{\tilde{K}_8\mathbf{x}^{(\tilde{y}_K)_6} - \tilde{K}_6 \mathbf{x}^{(\tilde{y}_K)_8} \Rightarrow k_{9} x_{13} - k_{18} x_{19}}
\end{array} \hspace{0.1in}
\begin{array}{c}
\underline{(3b) \& (4a)}: \vspace{0.05in} \\
\displaystyle{\tilde{K}_9\mathbf{x}^{(\tilde{y}_K)_{10}}- \tilde{K}_{10} \mathbf{x}^{(\tilde{y}_K)_9} \Rightarrow (k_{11} + k_{12}) x_{14} - k_{10} x_2x_4}\\
\displaystyle{\tilde{K}_{11}\mathbf{x}^{(\tilde{y}_K)_{12}} - \tilde{K}_{12} \mathbf{x}^{(\tilde{y}_K)_{11}} \Rightarrow (k_{14} + k_{15}) x_{18} - k_{13} x_5x_{11}}\\
\displaystyle{\tilde{K}_{12}\mathbf{x}^{(\tilde{y}_K)_{10}} - \tilde{K}_{10} \mathbf{x}^{(\tilde{y}_K)_{12}} \Rightarrow k_{12} x_{14} - k_{14} x_{18}}
\end{array} \vspace{0.1in}\\
\begin{array}{c}
\underline{(5a) \& (6b)}: \vspace{0.05in} \\
\displaystyle{\tilde{K}_{13}\mathbf{x}^{(\tilde{y}_K)_{14}} - \tilde{K}_{14} \mathbf{x}^{(\tilde{y}_K)_{13}} \Rightarrow (k_{20} + k_{21}) x_{15} - k_{19} x_5x_6}\\
\displaystyle{\tilde{K}_{15}\mathbf{x}^{(\tilde{y}_K)_{16}} - \tilde{K}_{16} \mathbf{x}^{(\tilde{y}_K)_{15}} \Rightarrow (k_{29} + k_{30}) x_{21} - k_{28} x_7x_{11}}\\
\displaystyle{\tilde{K}_{16}\mathbf{x}^{(\tilde{y}_K)_{14}} - \tilde{K}_{14} \mathbf{x}^{(\tilde{y}_K)_{16}} \Rightarrow k_{21} x_{15} - k_{30} x_{21}}
\end{array} \hspace{0.1in}
\begin{array}{c}
\underline{(5b) \& (6a)}: \vspace{0.05in} \\
\displaystyle{\tilde{K}_{17}\mathbf{x}^{(\tilde{y}_K)_{18}} - \tilde{K}_{18} \mathbf{x}^{(\tilde{y}_K)_{17}} \Rightarrow (k_{23} + k_{24}) x_{16} - k_{22} x_5x_7}\\
\displaystyle{\tilde{K}_{19}\mathbf{x}^{(\tilde{y}_K)_{20}} - \tilde{K}_{20} \mathbf{x}^{(\tilde{y}_K)_{19}} \Rightarrow (k_{26} + k_{27}) x_{20} - k_{25} x_8x_{11}}\\
\displaystyle{\tilde{K}_{20}\mathbf{x}^{(\tilde{y}_K)_{20}} - \tilde{K}_{18} \mathbf{x}^{(\tilde{y}_K)_{20}} \Rightarrow k_{24} x_{16} - k_{27} x_{20}}
\end{array}
\end{array}\]
\normalsize
As in Section \ref{section1}, the elements of the basis are grouped here to emphasize the connection with the linkage classes of the translation (\ref{translation2}).

\subsection{Application III: Network with negative feedback}
\label{section3}

Consider the following network:
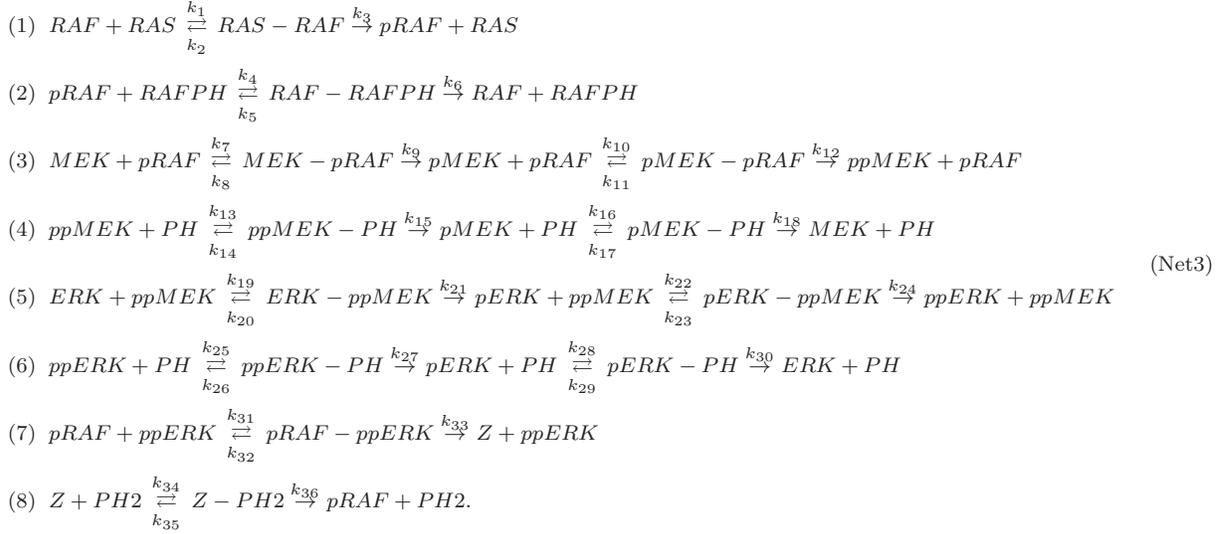
\scriptsize
\begin{equation}
\label{network3}\tag{Net3}
\begin{split}
(1) \; & \; RAF + RAS \; \mathop{\stackrel{k_1}{\rightleftarrows}}_{k_2} \; RAS-RAF \stackrel{k_3}{\rightarrow} pRAF + RAS\\
(2) \; & \; pRAF+RAFPH \; \mathop{\stackrel{k_4}{\rightleftarrows}}_{k_5} \; RAF-RAFPH \stackrel{k_6}{\rightarrow} RAF+RAFPH\\
(3) \; & \; MEK + pRAF \; \mathop{\stackrel{k_7}{\rightleftarrows}}_{k_8} \; MEK-pRAF \stackrel{k_9}{\rightarrow} pMEK+pRAF \; \mathop{\stackrel{k_{10}}{\rightleftarrows}}_{k_{11}} \; pMEK-pRAF \stackrel{k_{12}}{\rightarrow} ppMEK+pRAF\\
(4) \; & \; ppMEK+PH \; \mathop{\stackrel{k_{13}}{\rightleftarrows}}_{k_{14}} \; ppMEK-PH \stackrel{k_{15}}{\rightarrow} pMEK+PH \; \mathop{\stackrel{k_{16}}{\rightleftarrows}}_{k_{17}} \; pMEK-PH \stackrel{k_{18}}{\rightarrow} MEK + PH\\
(5) \; & \; ERK+ppMEK \; \mathop{\stackrel{k_{19}}{\rightleftarrows}}_{k_{20}} \; ERK-ppMEK \stackrel{k_{21}}{\rightarrow} pERK+ppMEK \; \mathop{\stackrel{k_{22}}{\rightleftarrows}}_{k_{23}} \; pERK-ppMEK \stackrel{k_{24}}{\rightarrow} ppERK+ppMEK\\
(6) \; & \; ppERK+PH \; \mathop{\stackrel{k_{25}}{\rightleftarrows}}_{k_{26}} \; ppERK-PH \stackrel{k_{27}}{\rightarrow} pERK+PH \; \mathop{\stackrel{k_{28}}{\rightleftarrows}}_{k_{29}} \; pERK-PH \stackrel{k_{30}}{\rightarrow} ERK+PH\\
(7) \; & \; pRAF + ppERK \; \mathop{\stackrel{k_{31}}{\rightleftarrows}}_{k_{32}} \; pRAF-ppERK \stackrel{k_{33}}{\rightarrow} Z + ppERK\\
(8) \; & \; Z + PH2 \; \mathop{\stackrel{k_{34}}{\rightleftarrows}}_{k_{35}} \; Z-PH2 \stackrel{k_{36}}{\rightarrow} pRAF + PH2.
\end{split}
\end{equation}
\normalsize
This is the network considered in Section \ref{section2} with the additional pathways (7) and (8). These pathways correspond to negative feedback since the fully phosphorylated extracellular signal-regulated kinase ($ppERK$) inhibits the activity of phosphorylated rapidly accelerated fibrosarcoma ($pRAF$) by binding it into the inactive form $Z$. The mechanism (\ref{network3}) has been considered by Asthagiri and Lauffenburger \cite{Asthagiri2001}, Dougherty \emph{et al.} \cite{Dougherty2005}, and Fritsche-Guenther \emph{et al.} \cite{fritsche2011}.

We now apply network translation to this network. We start by indexing the reactants as in P\'{e}rez Mill\'{a}n and Turjanski \cite{P-T}:

\footnotesize
\begin{equation}
\label{species3}
\begin{split}
& x_1 = [RAF], x_2 = [pRAF], x_3 = [MEK], x_4 = [pMEK], x_5 = [ppMEK], x_6 = [ERK], x_7 = [pERK],\\
& x_8 = [ppERK], x_9 = [RAS], x_{10} = [RAF P H], x_{11} = [P H], x_{12} = [PH2], x_{13} = [RAS - RAF],\\
& x_{14} = [MEK - pRAF], x_{15} = [pMEK - pRAF], x_{16} = [ERK - ppMEK], x_{17} = [pERK - ppMEK],\\
& x_{18} = [RAF - RAF PH], x_{19} = [ppMEK - P H], x_{20} = [pMEK - P H], x_{21} = [ppERK - P H],\\
& x_{22} = [pERK - P H],  x_{23} = [pRAF - ppERK], x_{24} = [Z - P H2], x_{25} = [Z].
\end{split}
\end{equation}
\normalsize
We now apply the following translation scheme where the reactants have been indexed according to (\ref{species3}):

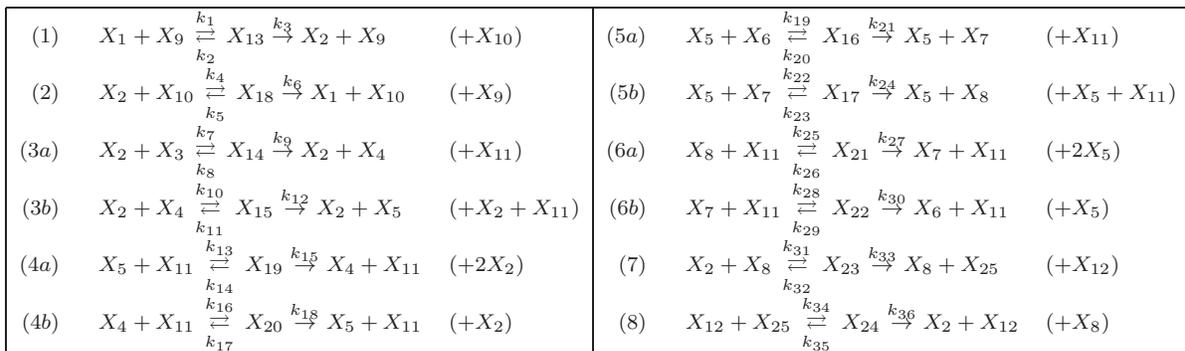
\scriptsize
\begin{equation}
\label{scheme3}\tag{Sch3}
\begin{array}{|rll|rll|}
\hline
(1) \; & \; \displaystyle{X_1 + X_9 \; \mathop{\stackrel{k_1}{\rightleftarrows}}_{k_2} \; X_{13} \stackrel{k_3}{\rightarrow} X_2 + X_9} & (+X_{10}) & 
(5a) \; & \; \displaystyle{X_5+X_6 \; \mathop{\stackrel{k_{19}}{\rightleftarrows}}_{k_{20}} \; X_{16} \stackrel{k_{21}}{\rightarrow} X_5+X_7} & (+X_{11})\\
(2) \; & \; \displaystyle{X_2+X_{10} \; \mathop{\stackrel{k_4}{\rightleftarrows}}_{k_5} \; X_{18} \stackrel{k_6}{\rightarrow} X_1+X_{10}} & (+X_{9}) &
(5b) \; & \; \displaystyle{X_5+X_7 \; \mathop{\stackrel{k_{22}}{\rightleftarrows}}_{k_{23}} \; X_{17} \stackrel{k_{24}}{\rightarrow} X_5+X_8} & (+X_5 + X_{11})\\
(3a) \; & \; \displaystyle{X_2 + X_3 \; \mathop{\stackrel{k_7}{\rightleftarrows}}_{k_8} \; X_{14} \stackrel{k_9}{\rightarrow} X_2+X_4} & (+X_{11}) &
(6a) \; & \; \displaystyle{X_8+X_{11} \; \mathop{\stackrel{k_{25}}{\rightleftarrows}}_{k_{26}} \; X_{21} \stackrel{k_{27}}{\rightarrow} X_7+X_{11}} & (+2X_5)\\
(3b) \; & \; \displaystyle{X_2+X_4 \; \mathop{\stackrel{k_{10}}{\rightleftarrows}}_{k_{11}} \; X_{15} \stackrel{k_{12}}{\rightarrow} X_2+X_5} & (+X_{2}+X_{11}) &
(6b) \; & \; \displaystyle{X_7+X_{11} \; \mathop{\stackrel{k_{28}}{\rightleftarrows}}_{k_{29}} \; X_{22} \stackrel{k_{30}}{\rightarrow} X_6+X_{11}} & (+X_5)\\
(4a) \; & \; \displaystyle{X_5+X_{11} \; \mathop{\stackrel{k_{13}}{\rightleftarrows}}_{k_{14}} \; X_{19} \stackrel{k_{15}}{\rightarrow} X_4+X_{11}} & (+2X_{2})&
(7) \; & \; \displaystyle{X_2+X_8 \; \mathop{\stackrel{k_{31}}{\rightleftarrows}}_{k_{32}} \; X_{23} \stackrel{k_{33}}{\rightarrow} X_8+X_{25}} & (+X_{12})\\
(4b) \; & \; \displaystyle{X_4+X_{11} \; \mathop{\stackrel{k_{16}}{\rightleftarrows}}_{k_{17}} \; X_{20} \stackrel{k_{18}}{\rightarrow} X_5 + X_{11}} & (+X_{2})&
(8) \; & \; \displaystyle{X_{12}+X_{25} \; \mathop{\stackrel{k_{34}}{\rightleftarrows}}_{k_{35}} \; X_{24} \stackrel{k_{36}}{\rightarrow} X_2+X_{12}} & (+X_8)\\
\hline
\end{array}
\end{equation}
\normalsize
The translation process (\ref{scheme3}) yields the following translated chemical reaction network $\tilde{\mathcal{N}} = (\tilde{\mathcal{S}},\tilde{\mathcal{C}},\tilde{\mathcal{C}}_K,\tilde{\mathcal{R}})$:

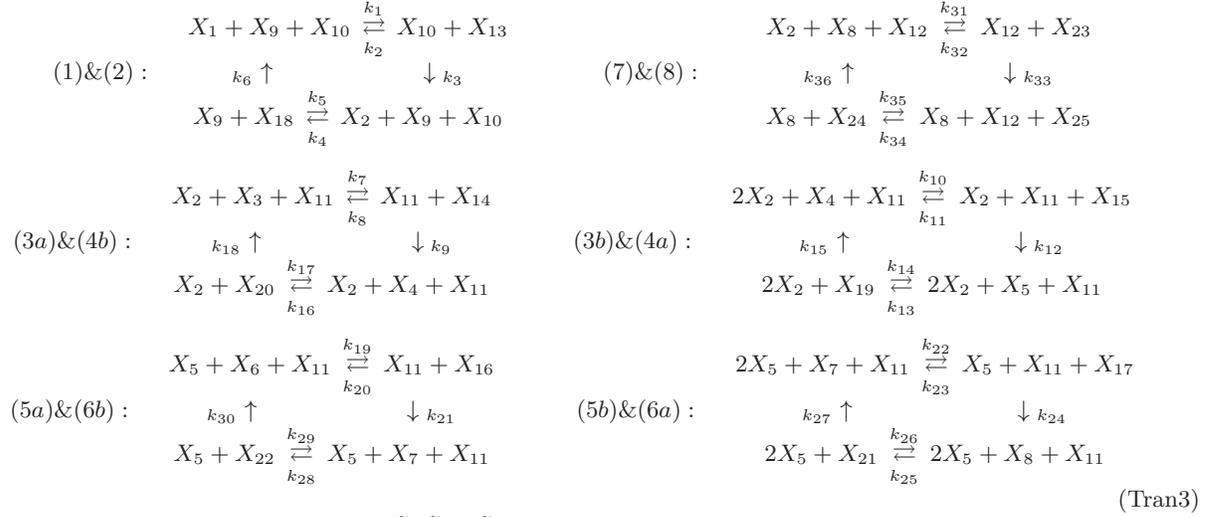
\footnotesize
\begin{equation}
\label{translation3}\tag{Tran3}
\begin{array}{c}
(1) \& (2): \; \; \; \begin{array}{c} \displaystyle{X_1 + X_9 + X_{10} \; \mathop{\stackrel{k_1}{\rightleftarrows}}_{k_2} \; X_{10} + X_{13}} \\ \displaystyle{{}_{k_6} \uparrow \hspace{2cm} \downarrow {}_{k_3}} \\ \displaystyle{X_9 + X_{18} \; \mathop{\stackrel{k_5}{\rightleftarrows}}_{k_4} \; X_2 + X_9 + X_{10}} \end{array} \; \; \; \; \; \; \; \; \; \; \; \;
(7) \& (8): \; \; \; \; \; \; \; \begin{array}{c} \displaystyle{X_2 + X_8 + X_{12} \; \mathop{\stackrel{k_{31}}{\rightleftarrows}}_{k_{32}} \; X_{12} + X_{23}} \\ \displaystyle{{}_{k_{36}} \uparrow \hspace{2cm} \downarrow {}_{k_{33}}} \\ \displaystyle{X_8 + X_{24} \; \mathop{\stackrel{k_{35}}{\rightleftarrows}}_{k_{34}} \; X_8 + X_{12} + X_{25}} \end{array}
 \vspace{0.1in} \\
(3a) \& (4b): \; \; \; \begin{array}{c} \displaystyle{X_2 + X_3 + X_{11} \; \mathop{\stackrel{k_7}{\rightleftarrows}}_{k_8} \; X_{11} + X_{14}} \\ \displaystyle{{}_{k_{18}} \uparrow \hspace{2cm} \downarrow {}_{k_{9}}} \\ \displaystyle{X_2 + X_{20} \; \mathop{\stackrel{k_{17}}{\rightleftarrows}}_{k_{16}} \; X_2 + X_4 + X_{11}} \end{array} \; \; \; \; \; \; \; \; \; \;
(3b) \& (4a): \; \; \; \begin{array}{c} \displaystyle{2X_2 + X_4 + X_{11} \; \mathop{\stackrel{k_{10}}{\rightleftarrows}}_{k_{11}} \; X_2 + X_{11} + X_{15}} \\ \displaystyle{{}_{k_{15}} \uparrow \hspace{2.2cm} \downarrow {}_{k_{12}}} \\ \displaystyle{2X_2 + X_{19} \; \mathop{\stackrel{k_{14}}{\rightleftarrows}}_{k_{13}} \; 2X_2 + X_5 + X_{11}} \end{array} \vspace{0.1in} \\
(5a) \& (6b): \; \; \; \begin{array}{c} \displaystyle{X_5 + X_6 + X_{11} \; \mathop{\stackrel{k_{19}}{\rightleftarrows}}_{k_{20}} \; X_{11} + X_{16}} \\ \displaystyle{{}_{k_{30}} \uparrow \hspace{2cm} \downarrow {}_{k_{21}}} \\ \displaystyle{X_5 + X_{22} \; \mathop{\stackrel{k_{29}}{\rightleftarrows}}_{k_{28}} \; X_5 + X_7 + X_{11}} \end{array} \; \; \; \; \; \; \; \; \; \;
(5b) \& (6a): \; \; \; \begin{array}{c} \displaystyle{2X_5 + X_7 + X_{11} \; \mathop{\stackrel{k_{22}}{\rightleftarrows}}_{k_{23}} \; X_5 + X_{11} + X_{17}} \\ \displaystyle{{}_{k_{27}} \uparrow \hspace{2.2cm} \downarrow {}_{k_{24}}} \\ \displaystyle{2X_5 + X_{21} \; \mathop{\stackrel{k_{26}}{\rightleftarrows}}_{k_{25}} \; 2X_5 + X_8 + X_{11}} \end{array}
\end{array}
\end{equation}
\normalsize
and the associated kinetic graph $(\tilde{\mathcal{S}},\tilde{\mathcal{C}}_K,\tilde{\mathcal{R}})$:

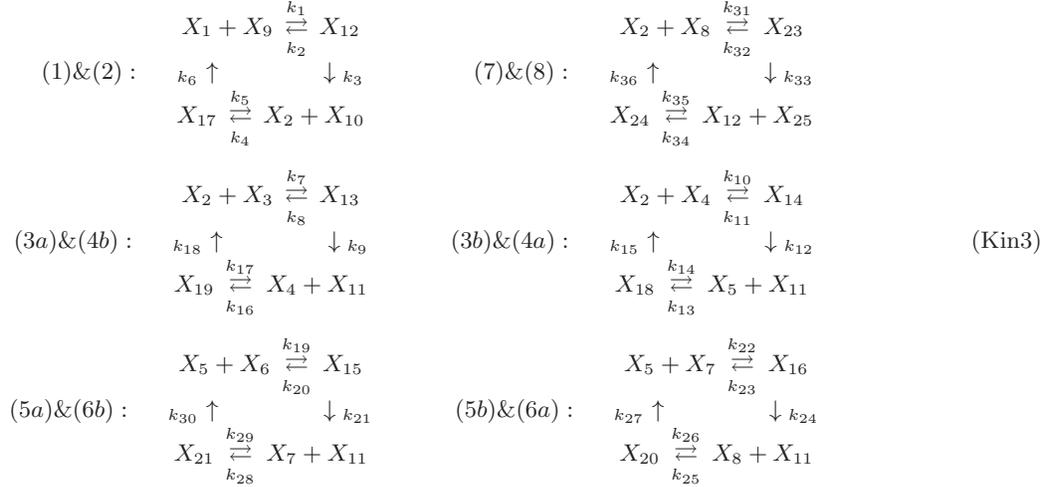
\footnotesize
\begin{equation}
\label{kinetic3}\tag{Kin3}
\begin{array}{c}
\; \; \; \; (1) \& (2): \; \; \; \begin{array}{c} \displaystyle{X_1 + X_9 \; \mathop{\stackrel{k_1}{\rightleftarrows}}_{k_2} \; X_{12}} \\ \displaystyle{{}_{k_6} \uparrow \hspace{1.4cm} \downarrow {}_{k_3}} \\ \displaystyle{X_{17} \; \mathop{\stackrel{k_5}{\rightleftarrows}}_{k_4} \; X_2 + X_{10}} \end{array} \; \; \; \; \; \; \; \; \; \;
\; \; \; \; (7) \& (8): \; \; \;\begin{array}{c} \displaystyle{X_2 + X_8 \; \mathop{\stackrel{k_{31}}{\rightleftarrows}}_{k_{32}} \; X_{23}} \\ \displaystyle{{}_{k_{36}} \uparrow \hspace{1.4cm} \downarrow {}_{k_{33}}} \\ \displaystyle{X_{24} \; \mathop{\stackrel{k_{35}}{\rightleftarrows}}_{k_{34}} \; X_{12} + X_{25}} \end{array} \vspace{0.1in} \\
(3a) \& (4b): \; \; \; \begin{array}{c} \displaystyle{X_2 + X_3 \; \mathop{\stackrel{k_7}{\rightleftarrows}}_{k_8} \; X_{13}} \\ \displaystyle{{}_{k_{18}} \uparrow \hspace{1.4cm} \downarrow {}_{k_{9}}} \\ \displaystyle{X_{19} \; \mathop{\stackrel{k_{17}}{\rightleftarrows}}_{k_{16}} \; X_4 + X_{11}} \end{array} \; \; \; \; \; \; \; \; \; \;
(3b) \& (4a): \; \; \; \begin{array}{c} \displaystyle{X_2 + X_4 \; \mathop{\stackrel{k_{10}}{\rightleftarrows}}_{k_{11}} \; X_{14}} \\ \displaystyle{{}_{k_{15}} \uparrow \hspace{1.4cm} \downarrow {}_{k_{12}}} \\ \displaystyle{X_{18} \; \mathop{\stackrel{k_{14}}{\rightleftarrows}}_{k_{13}} \; X_5 + X_{11}} \end{array} \vspace{0.1in} \\
(5a) \& (6b): \; \; \; \begin{array}{c} \displaystyle{X_5 + X_6 \; \mathop{\stackrel{k_{19}}{\rightleftarrows}}_{k_{20}} \; X_{15}} \\ \displaystyle{{}_{k_{30}} \uparrow \hspace{1.4cm} \downarrow {}_{k_{21}}} \\ \displaystyle{X_{21} \; \mathop{\stackrel{k_{29}}{\rightleftarrows}}_{k_{28}} \; X_7 + X_{11}} \end{array} \; \; \; \; \; \; \; \; \; \;
(5b) \& (6a): \; \; \; \begin{array}{c} \displaystyle{X_5 + X_7 \; \mathop{\stackrel{k_{22}}{\rightleftarrows}}_{k_{23}} \; X_{16}} \\ \displaystyle{{}_{k_{27}} \uparrow \hspace{1.4cm} \downarrow {}_{k_{24}}} \\ \displaystyle{X_{20} \; \mathop{\stackrel{k_{26}}{\rightleftarrows}}_{k_{25}} \; X_8 + X_{11}} \end{array}
\end{array}
\end{equation}
\normalsize
As for the networks in Sections \ref{section1} and \ref{section2}, we can determine that the translation $\tilde{\mathcal{N}}$ is proper, weakly reversible, and has both $\tilde{\delta} = 0$ and $\tilde{\delta}_K = 0$. We make the following enumeration of the kinetic complexes:

\footnotesize
\begin{equation}
\label{kineticcomplexes2}
\begin{split}
(1) \& (2): \; \; \; \; \; & (\tilde{C}_K)_1 = X_1 + X_9, (\tilde{C}_K)_2 = X_{13}, (\tilde{C}_K)_3 = X_2+X_{10}, (\tilde{C}_K)_4 = X_{18},\\
(3a) \& (4b): \; \; \; \; \; & (\tilde{C}_K)_5 = X_2 + X_3, (\tilde{C}_K)_6 = X_{14}, (\tilde{C}_K)_7 = X_4 + X_{11}, (\tilde{C}_K)_8 = X_{20},\\
(3b) \& (4a): \; \; \; \; \; & (\tilde{C}_K)_9 = X_2 + X_4, (\tilde{C}_K)_{10} = X_{15}, (\tilde{C}_K)_{11} = X_5 + X_{11}, (\tilde{C}_K)_{12} = X_{19},\\
(5a) \& (6b): \; \; \; \; \; & (\tilde{C}_K)_{13} = X_5 + X_6, (\tilde{C}_K)_{14} = X_{16}, (\tilde{C}_K)_{15} = X_7 + X_{11}, (\tilde{C}_K)_{16} = X_{22}\\
(5b) \& (6a): \; \; \; \; \; & (\tilde{C}_K)_{17} = X_5 + X_7, (\tilde{C}_K)_{18} = X_{18}, (\tilde{C}_K)_{19} = X_8 + X_{11}, (\tilde{C}_K)_{20} = X_{21},\\
(7) \& (8): \; \; \; \; \; & (\tilde{C}_K)_{21} = X_2 + X_8, (\tilde{C}_K)_{22} = X_{23}, (\tilde{C}_K)_{23} = X_{12}+X_{25}, (\tilde{C}_K)_{24} = X_{24},
\end{split}
\end{equation}
\normalsize
and compute the following tree constants:

\footnotesize
\begin{equation}
\begin{array}{|l|l|l|l|}
\hline
\tilde{K}_1 = (k_2 + k_3)k_4k_6 & \tilde{K}_7 = k_7k_9(k_{17}+k_{18}) & \tilde{K}_{13} = (k_{20}+k_{21})k_{28}k_{30} & \tilde{K}_{19} = k_{22}k_{24}(k_{26}+k_{27})\\
\tilde{K}_2 = k_1k_4k_6 & \tilde{K}_8 = k_7k_9k_{16} & \tilde{K}_{14} = k_{19}k_{28}k_{30} & \tilde{K}_{20} = k_{23}k_{24}k_{25}\\
\tilde{K}_3 = k_1k_3(k_5+k_6) & \tilde{K}_9 = (k_{11}+k_{12})k_{13}k_{15} & \tilde{K}_{15} = k_{19}k_{21}(k_{29}+k_{30}) & \tilde{K}_{21} = (k_{32}+k_{33})k_{34}k_{36}\\
\tilde{K}_4 = k_1k_3k_4 & \tilde{K}_{10} = k_{10}k_{13}k_{15} & \tilde{K}_{16} = k_{19}k_{21}k_{28} & \tilde{K}_{22} = k_{31}k_{34}k_{36}\\
\tilde{K}_5 = (k_8+k_9)k_{16}k_{18} & \tilde{K}_{11} = k_{10}k_{12}(k_{14}+k_{15}) & \tilde{K}_{17} = (k_{23}+k_{24})k_{25}k_{27} & \tilde{K}_{23} = k_{31}k_{33}(k_{35}+k_{36})\\
\tilde{K}_6 = k_7k_{18}k_{16} &\tilde{K}_{12} = k_{10}k_{12}k_{13} & \tilde{K}_{18} = k_{22}k_{25}k_{27} & \tilde{K}_{24} = k_{31}k_{33}k_{34}.\\
\hline
\end{array}
\end{equation}
\normalsize
By property 2. of  Theorem 5 of Johnston \cite{J1}, a basis of the steady state ideal can then be given by:

\footnotesize
\[\begin{array}{c}
\begin{array}{c}
\underline{(1) \& (2)}: \vspace{0.05in} \\
\displaystyle{\tilde{K}_1\mathbf{x}^{(\tilde{y}_K)_2} - \tilde{K}_2 \mathbf{x}^{(\tilde{y}_K)_1} \Rightarrow (k_2 + k_3) x_{13} - k_1 x_1x_9}\\
\displaystyle{\tilde{K}_3\mathbf{x}^{(\tilde{y}_K)_4} - \tilde{K}_4 \mathbf{x}^{(\tilde{y}_K)_3} \Rightarrow (k_5 + k_6) x_{18} - k_4 x_2x_{10}}\\
\displaystyle{\tilde{K}_4\mathbf{x}^{(\tilde{y}_K)_2} - \tilde{K}_2 \mathbf{x}^{(\tilde{y}_K)_4} \Rightarrow k_3 x_{13} - k_6 x_{18}}
\end{array} \hspace{0.1in}
\begin{array}{c}
\underline{(7) \& (8)}: \vspace{0.05in} \\
\displaystyle{\tilde{K}_{21}\mathbf{x}^{(\tilde{y}_K)_{22}} - \tilde{K}_{22} \mathbf{x}^{(\tilde{y}_K)_{21}} \Rightarrow (k_{32} + k_{33}) x_{23} - k_{31} x_2x_8}\\
\displaystyle{\tilde{K}_{23}\mathbf{x}^{(\tilde{y}_K)_{24}} - \tilde{K}_{24} \mathbf{x}^{(\tilde{y}_K)_{23}} \Rightarrow (k_{35} + k_{36}) x_{24} - k_{34} x_{12}x_{25}}\\
\displaystyle{\tilde{K}_{24}\mathbf{x}^{(\tilde{y}_K)_{22}} - \tilde{K}_{22} \mathbf{x}^{(\tilde{y}_K)_{24}} \Rightarrow k_{33} x_{23} - k_{36} x_{24}}
\end{array} \vspace{0.1in}\\
\begin{array}{c}
\underline{(3a) \& (4b)}: \vspace{0.05in} \\
\displaystyle{\tilde{K}_5\mathbf{x}^{(\tilde{y}_K)_6} - \tilde{K}_6 \mathbf{x}^{(\tilde{y}_K)_5} \Rightarrow (k_8 + k_9) x_{14} - k_7 x_2x_3}\\
\displaystyle{\tilde{K}_7\mathbf{x}^{(\tilde{y}_K)_8} - \tilde{K}_8 \mathbf{x}^{(\tilde{y}_K)_7} \Rightarrow (k_{17} + k_{18}) x_{20} - k_{16} x_4x_{11}}\\
\displaystyle{\tilde{K}_8\mathbf{x}^{(\tilde{y}_K)_6} - \tilde{K}_6 \mathbf{x}^{(\tilde{y}_K)_8} \Rightarrow k_{9} x_{14} - k_{18} x_{20}}
\end{array} \hspace{0.1in}
\begin{array}{c}
\underline{(3b) \& (4a)}: \vspace{0.05in} \\
\displaystyle{\tilde{K}_9\mathbf{x}^{(\tilde{y}_K)_{10}}- \tilde{K}_{10} \mathbf{x}^{(\tilde{y}_K)_9} \Rightarrow (k_{11} + k_{12}) x_{15} - k_{10} x_2x_4}\\
\displaystyle{\tilde{K}_{11}\mathbf{x}^{(\tilde{y}_K)_{12}} - \tilde{K}_{12} \mathbf{x}^{(\tilde{y}_K)_{11}} \Rightarrow (k_{14} + k_{15}) x_{19} - k_{13} x_5x_{11}}\\
\displaystyle{\tilde{K}_{12}\mathbf{x}^{(\tilde{y}_K)_{10}} - \tilde{K}_{10} \mathbf{x}^{(\tilde{y}_K)_{12}} \Rightarrow k_{12} x_{15} - k_{15} x_{19}}
\end{array} \vspace{0.1in}\\
\begin{array}{c}
\underline{(5a) \& (6b)}: \vspace{0.05in} \\
\displaystyle{\tilde{K}_{13}\mathbf{x}^{(\tilde{y}_K)_{14}} - \tilde{K}_{14} \mathbf{x}^{(\tilde{y}_K)_{13}} \Rightarrow (k_{20} + k_{21}) x_{16} - k_{19} x_5x_6}\\
\displaystyle{\tilde{K}_{15}\mathbf{x}^{(\tilde{y}_K)_{16}} - \tilde{K}_{16} \mathbf{x}^{(\tilde{y}_K)_{15}} \Rightarrow (k_{29} + k_{30}) x_{22} - k_{28} x_7x_{11}}\\
\displaystyle{\tilde{K}_{16}\mathbf{x}^{(\tilde{y}_K)_{14}} - \tilde{K}_{14} \mathbf{x}^{(\tilde{y}_K)_{16}} \Rightarrow k_{21} x_{16} - k_{30} x_{22}}
\end{array} \hspace{0.1in}
\begin{array}{c}
\underline{(5b) \& (6a)}: \vspace{0.05in} \\
\displaystyle{\tilde{K}_{17}\mathbf{x}^{(\tilde{y}_K)_{18}} - \tilde{K}_{18} \mathbf{x}^{(\tilde{y}_K)_{17}} \Rightarrow (k_{23} + k_{24}) x_{17} - k_{22} x_5x_7}\\
\displaystyle{\tilde{K}_{19}\mathbf{x}^{(\tilde{y}_K)_{20}} - \tilde{K}_{20} \mathbf{x}^{(\tilde{y}_K)_{19}} \Rightarrow (k_{26} + k_{27}) x_{21} - k_{25} x_8x_{11}}\\
\displaystyle{\tilde{K}_{20}\mathbf{x}^{(\tilde{y}_K)_{20}} - \tilde{K}_{18} \mathbf{x}^{(\tilde{y}_K)_{20}} \Rightarrow k_{24} x_{17} - k_{27} x_{21}}
\end{array}
\end{array}\]
\normalsize
As in Sections \ref{section1} and \ref{section2}, the elements of the basis are grouped here to emphasize the connection with the linkage classes of the translation (\ref{translation3}).

\section{Conclusions}

In this note, we have illustrated how network translation (Johnston \cite{J1}) can be used to characterize mass action systems with toric steady states (P\'{e}rez Mill\'{a}n \emph{et al.} \cite{M-D-S-C}) by applying it to three related MAPK networks known to have toric steady states (P\'{e}rez Mill\'{a}n and Turjanski \cite{P-T}). The process of network translation involves corresponding the original network to a generalized network (M\"{u}ller and Regensburger \cite{M-R}) with superior structure. In all three cases considered here, we were able to construct a translation which was proper, weakly reversible, and for which $\tilde{\delta} = 0$ and $\tilde{\delta}_K = 0$. This is sufficient to guarantee ``complex-balanced steady states'' for the translated system and to guarantee toric steady states of the original system.

We make a few notes regarding the scope of this work and avenues for future work:
\begin{enumerate}
\item
Network translation is not the only approach capable of guaranteeing the existence of toric steady states. Such states may be guaranteed by computing the kernel of a particular structural matrix, as illustrated by P\'{e}rez Mill\'{a}n \emph{et al.} \cite{M-D-S-C}, or through direct computation of the Gr\"{o}bner basis of the mass action system (see Cox \emph{et al.} \cite{C-L-O}). The contribution made in this note, and by network translation in general, is in relating the algebraic details of the steady state analysis explicitly to the topological structure of the underlying chemical reaction network. This is in keeping with the flavor of so-called \emph{chemical reaction network theory} which has placed considerable emphasis of dynamical properties which follow from properties of the network structure itself (see Feinberg \cite{F1}, Horn \cite{H}, and Horn and Jackson \cite{H-J1}).
\item
The translation schemes (\ref{scheme1}), (\ref{scheme2}), and (\ref{scheme3}) were chosen due to the similarity in network structures with that of the network considered in Application II of Johnston \cite{J1}. In general practice, however, such a scheme may not be immediately apparent. Recent work has been conducted by Johnston \cite{J2} on computational algorithms for constructing network translations. The algorithm given there is capable of corresponding regular mass action systems to generalized ones which have the same steady states. To date, however, the algorithm requires a given set of potential stoichiometric translated complexes which may not be obvious. Research on how to algorithmically determine the network structure of translations is therefore ongoing.
\item
We re-iterate that this note only addresses the construction of the basis for the steady state ideal and does not address the multistationarity or multistability of the system. To date, no results are known which relate the process of network translation to a network's capacity for multistationarity. Given that conditions are known on both regular (Craciun and Feinberg \cite{C-F1,C-F2}; P\'{e}rez Mill\'{a}n \emph{et al.} \cite{M-D-S-C}; M\"{u}ller \emph{et al.} \cite{M-F-R-C-S-D}) and generalized mass action systems (M\"{u}ller and Regensberger \cite{M-R}) which guarantee multistationarity for some set of rate constants, we therefore ask the following: Are there properties of the network translation process itself which guarantee, or prohibit, that the underlying network has the capacity for multistationarity?
\end{enumerate}


\end{document}